\documentclass{amsart}
\usepackage{amsfonts}
\usepackage{amssymb}
\numberwithin{equation}{section}
\usepackage{amscd}

\def\longnearrow{\hbox{{\begin{Large}$\displaystyle\nearrow$\end{Large}}}}
\def\longsearrow{\hbox{{\begin{Large}$\displaystyle\searrow$\end{Large}}}}
\def\spacebox{\lower2ex\hbox{\vbox to 22pt{\vfill}}}
\def\littlespacebox{\lower2ex\hbox{\vbox to 16pt{\vfill}}}

\newcommand{\A}{{\mathcal A}}
\newcommand{\B}{{\mathcal B}}

\newcommand{\SL}{{\mathcal L}}

\newcommand{\N}{{\mathbb N}}

\newcommand{\OL}{\mathcal {OL}}
\newcommand{\al}{\alpha}
\newcommand{\be}{\beta}
\newcommand{\de}{\delta}
\newcommand{\e}{\varepsilon}
\newcommand{\la}{\lambda}
\newcommand{\ga}{\gamma}
\newcommand{\p}{\varphi}
\newcommand{\ten}{\check\otimes}
\makeatletter
\mathchardef\hugecheck="7014
\newcommand\hugesize{\@setfontsize\hugesize{25pt}{0}}
\newcommand\smallhugesize{\@setfontsize\smallhugesize{20pt}{0}}

\def\specialchecksmall{\mbox{\hbox to
0pt{\raisebox{-4pt}{\smallhugesize$\hugecheck$}}
                                $\kern-2.7pt\otimes$}}

\makeatother

\newtheorem{claim}{}[section]
\newtheorem{theorem}[claim]{Theorem}
\newtheorem{lemma}[claim]{Lemma}
\newtheorem{proposition}[claim]{Proposition}

\newtheorem{question}[claim]{Question}
\newtheorem{conjecture}[claim]{Conjecture}

\def\proclaim #1. #2\par{\medbreak
\noindent{\bf#1.\enspace}{\sl#2}\par\medbreak}
\begin{document}

\title[On ${\OL}_{\infty}$ structure of nuclear $C^*$-algebras]
{On ${\OL}_{\infty}$ structure of nuclear $C^*$-algebras}
\date{July 14, 2000, revised June 2,  2001}
\subjclass{Primary 46L07 and 46L05; Secondary 47L25}
\author{Marius Junge}
\address{Department of Mathematics\\
University of Illinois, Urbana, IL 61801 USA}
\email[Marius Junge]{junge@math.uiuc.edu}
\author{Narutaka Ozawa}
\address{Department of Mathematics\\
Texas A\&M University, College Station, TX 77843 USA\\
and
Department of Mathematical Sciences\\
University of Tokyo, Komaba, 153-8914, Japan}
\email{ozawa@math.tamu.edu}
\author{Zhong-Jin Ruan}
\address{Department of Mathematics\\
University of Illinois, Urbana, IL 61801 USA}
\email[Zhong-Jin Ruan]{ruan@math.uiuc.edu}
\thanks{Junge and Ruan were partially supported by the National Science
Foundation}
\thanks{Ozawa was supported by the Japanese Society for Promotion of
Science}

\begin{abstract}
We study the local operator space structure of nuclear $C^*$-algebras.
It is shown that a $C^*$-algebra  is nuclear if and only if
it is an $\OL_{\infty, \la}$ space for some (and actually for every)
$\la > 6$.
The $\OL_\infty$ constant $\lambda$ provides an interesting
invariant
\[
\OL_\infty (\A) = \inf\{ \la: ~ \A ~{\rm is ~ an} ~ {\OL}_{\infty, \la}
~ {\rm space }\}
\]
for nuclear $C^*$-algebras.
Indeed, if $\A$ is a nuclear $C^*$-algebra, then we have
$1\le \OL_\infty (\A)  \le 6$,
and  if $\A$ is a unital nuclear $C^*$-algebra with
$\OL_{\infty} (\A) \le  (\frac {1+{\sqrt 5}}2)^{\frac 12}$,
we show that  $\A$ must be stably  finite.

We also investigate the connection between the rigid $\OL_{\infty, 1^+}$
structure and  the rigid complete order $\OL_{\infty, 1^+}$  structure on
$C^*$-algebras, where the  latter structure   has been studied by
Blackadar and Kirchberg in their characterization of strong NF $C^*$-algebras.
Another main result of this paper is to show that these two local
structrues are  actually equivalent  on  unital nuclear $C^*$-algebras.
We obtain this by showing that if a  unital (nuclear) $C^*$-algebra is a rigid
${\OL}_{\infty, 1^+}$  space,  then it is inner quasi-diagonal, and thus is  a
strong NF algebra.
It is also shown that if  a unital (nuclear) $C^*$-algebra is an
${\OL}_{\infty, 1^+}$ space,  then it is  quasi-diagonal, and thus is an NF
algebra.
\end{abstract}

\maketitle
\let\text=\mbox
\let\cal=\mathcal
%\linespread{1.4}

\section{Introduction}

The purpose of this paper is to  use operator space theory to
study the  underlying  operator space structure (more precisely, the
${\OL}_{\infty}$  structure) of  nuclear $C^*$-algebras.
Roughly speaking, an ${\OL}_{\infty}$ space is an operator space which
can be ``locally paved up'' by finite-dimensional  $C^*$-algebras.
This notion is a natural operator space analogue of ${\SL}_{\infty}$
spaces introduced by  Lindenstrauss  and Pe{\l}czy\'{n}ski \cite{LP}.
As in Banach space theory, there are two typical ways to pave up
an operator space by finite-dimensional $C^*$-algebras, i.e.
by uniformly completely bounded isomorphic injections, or by
completely isometric (i.e. rigid) injections.
We will see that it is crucially important to distinguish these different
local structures for nuclear $C^*$-algebras.

To explain our motivation, let  us first recall some Banach space results from
\cite{LP}.
Given $1\leq p\leq \infty$ and $\lambda >1,$
a Banach space $V$ is said to be an $\mathcal{L}_{p,\lambda }$
\emph{space}
if for every finite-dimensional  subspace $E$ of $V$,
there exists a finite-dimensional subspace $F$ of $V$ such that
$E\subseteq F$ (with $\mbox{dim} F = n$)
and  the Banach-Mazur distance
\[
d (F, \ell_p(n)) = \inf \left\{ \Vert T\Vert \Vert T^{-1}\Vert
:T:F\rightarrow \ell_p(n)
\text{ a linear isomorphism}\right\}< \la.
\]
A Banach  space $V$ is said to be an ${\SL}_{p, 1^+}$ \emph{space}  if
it is an ${\SL}_{p, \la }$ space for every $\la > 1$,
and is said to be a \emph {rigid} $\SL_{p, 1^+}$
\emph {space} if there exists a collection of finite-dimensional
subspaces $F_\al$ of $V$ such that each $F_\al$ is isometric to an
$\ell_p(n_\al)$ space and  the union of $F_\al$ is norm dense in $V$.

It is well-known (see \cite{LP}) that for $1\leq p<\infty$, a Banach
space $V$ is an $\mathcal{L}_{p,1^+}$ space  if and only if it is
isometric to an $L_p(X, {\mathcal M}, \mu)$ space (and thus is a rigid
$\mathcal{L}_{p,1^+}$ space).
Therefore, ${\SL}_{p, \la}$ spaces are natural local generalization of
$L_p(X, {\mathcal M}, \mu)$ spaces. The situation is
more subtle for $p = \infty$.
It is easy to see that commutative $C^*$-algebras $C(\Omega)$
(these include $L_\infty(X, {\mathcal M}, \mu)$ spaces) are
${\SL}_{\infty, 1^+}$ spaces.
However, there are many other $\mathcal{L}_{\infty ,1^+}$ spaces.
In general, it is known that a Banach space is an $\mathcal{L}_{\infty,1^+}$
space if and only if it is a predual of some $L_1(X, {\mathcal M}, \mu)$ space,
and it is a non-trivial result of  Michael and Pe{\l}czy\'{n}ski \cite {MP}
that a Banach space $V$ is an ${\SL}_{\infty, 1^+}$ space if and only if
it is a rigid ${\SL}_{\infty, 1^+}$ space.

Operator spaces are natural non-commutative quantization of Banach spaces.
An operator space can be (concretely) defined to be a norm closed linear space
of operators on some Hilbert space $H$, which is equipped
with a distinguished  \emph{matrix norm} obtained from $\B(H)$.
The morphisms between operator spaces are \emph{completely bounded mappings}.
There are many parallel results as well as many distinctions between
operator spaces and Banach spaces.
Nevertheless,  Banach space theory always provides an important
source of inspiration for  the development of operator space theory.
The readers are referred to the recent book of Effros and Ruan \cite {ERbook}
and the book of Pisier \cite {Pibook} for details.

The operator space analogue of ${\SL}_{p}$ spaces was first studied by
Effros and Ruan \cite {ER}.
Let us first  recall from Pisier  \cite {Pip} that if $\B$ is a
finite-dimensional $C^*$-algebra, then
we may use complex interpolation to obtain a natural operator
space structure on the non-commutative $L_p(\B)$ space
(for  $1\le p \le \infty$).
In particular, we have $L_\infty(\B) = \B$ and $L_1(\B) = \B_*$.
An operator space $V$ is said to be an ${\OL}_{p, \la}$ space
(for some $\la > 1$) if for every finite-dimensional  subspace
$E$ of $V$, there exists a finite-dimensional $C^*$-algebra $\B$
and a finite-dimensional subspace $F$ of
$V$ such that $E\subseteq F$  and the completely bounded Banach-Mazur
distance
\begin{equation}
\label {F1.1a}
d_{cb}(F, L_p(\B)) < \la,
\end{equation}
where (\ref {F1.1a}) is equivalent to saying that there exists
  a linear isomorphism
\[
\p :  L_p(\B) \to F
\]
such that $\|\p\|_{cb}  \|\p^{-1}\|_{cb} < \la.$
An operator space $V$ is said to be an ${\OL}_{p, 1^+}$ \emph{space}
if it is an ${\OL}_{p, \la}$  space for every $\la > 1$,
and is said to be
a \emph{rigid} ${\OL}_{p, 1^+}$  \emph{space} if
for every  $x_1,\cdots, x_n \in V$ and $\e > 0$, there exists a
finite-dimensional $C^*$-algebra $\B$ and
a  complete isometry $\p$ from $L_p(\B)$ onto a finite-dimensional
subspace $F$ of $V$ such that
\[
dist(x_i, F)  < \e
\]
for all $i = 1, \cdots , n$.
Using a standard perturbation argument, we can easily show that
rigid  ${\OL}_{p, 1^+}$ spaces are  ${\OL}_{p, 1^+}$ spaces.
For $p=\infty$,  this has been discussed  in \cite {ER}.

${\OL}_{1, \la}$ spaces have been intensively studied in
\cite {ER} and \cite {NO}.
It was shown in \cite {ER}  that if
$V = R_*$ is the  operator predual of a von Neumann algebra $R$
on a separable Hilbert space,
then  $V$ is an  ${\OL}_{1, \la}$ space for some $\la > 1$
(respectively, an ${\OL}_{1, 1^+}$ space or  a rigid ${\OL}_{1, 1^+}$ space)
if and only if $R$ is an  injective von Neumann algebra.
The separability can be removed by a result of Haagerup, which was
stated and proved in \cite [Appendix]{G4}.
This shows that various  notions of $\OL_1$ structures
are all equivalent  on the operator preduals of von Neumann algebras.
Recently, Ng and Ozawa \cite {NO} have  proved  that a separable
operator space $V$ is an
${\OL}_{1, 1^+}$ space if and only if
$V$  is completely isometric to the operator predual of an injective von
Neumann algebra.
Therefore, the local structure of the operator preduals of injective von
Neumann algebras has been completely understood.

Turning to the ${\OL}_{\infty, \la}$ space case,
it was observed in \cite {ER} that every ${\OL}_{\infty, \la}$ space
is $\la$-\emph{nuclear}, i.e.
there exist  approximate diagrams of completely bounded mappings
\begin{equation}
\label{F1.nuc}
\begin{array}{ccccc}
&  & M _{n(\gamma)} &  &  \\
& {\scriptstyle \al_\gamma }\longnearrow  &  & \longsearrow
{\scriptstyle
\be_\gamma} &  \\
V\hspace{-10 pt} &  & \stackrel{id_V}{\hspace{-10pt}{\hbox to 30pt
{\rightarrowfill}}\hspace{-10pt} } &  &\hspace{-10 pt} V
\end{array}
\end{equation}
such that $\|\al_\gamma\|_{cb} \|\be_\gamma\|_{cb} \le \la$ and
$\be_\gamma \circ \al_\gamma\to id_V$ in the point-norm topology.
If the mappings $\al_\ga$ and $\be_\ga$ in (\ref {F1.nuc}) are
completely contractive, we say that $V$ is a {\em nuclear} operator
space. The notion  of nuclear operator space  first appeared in
Kirchberg \cite {Ki2}. Smith \cite {Smith} showed that for
$C^*$-algebras, this is equivalent to the usual definition
introduced by  Lance \cite {Lance}. In fact, Pisier \cite {PiOH}
proved that for $C^*$-algebras, the nuclearity is equivalent to
the $\la$-nuclearity. Therefore,  if a $C^*$-algebra $\A$ is an
${\OL}_{\infty, \la}$ space  for some $\la > 1$, then $\A$ must be
nuclear (see \cite [Proposition 4.9]{ER}).
Surprisingly, the local structure of nuclear $C^*$-algebras turns out to
be more sophisticated.
The aim of this paper is to study various notions of $\OL_\infty$ structures
and related properties on nuclear $C^*$-algebras.

Let us first look at the rigid $\OL_{\infty, 1^+}$ structure on
$C^*$-algebras.
It is easy to see that if a $C^*$-algebra $\A$ is a rigid
$\OL_{\infty, 1^+}$ space, then each complete isometry
$\p_{F, \e}^{-1}: F \to \B$  (given in the definition)
extends to a complete  contraction $\psi_{F, \e} : \A \to \B$.
Then we may obtain a net of finite-rank completely contractive projections
\begin{equation}
P_{F, \e} = \p_{F, \e} \circ \psi_{F, \e}
\end{equation}
on $\A$, which converges to $id_\A$ in the point-norm topology.
We show in $\S 2$ (see Theorem \ref {PR.2}) that the existence of
such a net of finite-rank completely contractive projections on
a unital $C^*$-algebra $\A$ implies that $\A$ is a
rigid  $\OL_{\infty, 1^+}$ space.
We note that, in general, the range space
$P(\A)$ of a finite-rank completely contractive projection $P$ on $\A$ is
not necessarily a finite-dimensional $C^*$-algebra.
It is a finite-dimensional ternary ring of operators (see Youngson \cite {Yo}),
and thus is a finite direct sum of rectangular matrices
(see Smith \cite {Smith2}).
This motivated us to consider \emph{ rigid rectangular}
$\OL_{\infty, 1^+}$ \emph{spaces} in $\S 2$ and to show in
Theorem \ref {PR.2} that a unital $C^*$-algebra  is a
rigid $\OL_{\infty, 1^+}$ space if and only if it is a
rigid rectangular $\OL_{\infty, 1^+}$ space.
The theory of ternary ring of operators plays a key role in $\S 2$.

It is well-known that for every $C^*$-algebra $\A$, there is a canonical
\emph{matrix order} on $\A$ given by the positive cones  $M_n(\A)^+$
in the matrix spaces $M_n(\A)$ for $n \in {\Bbb N}$.
Then it is natural to consider the \emph{rigid complete order}
$\OL_{\infty, 1^+}$ structure on $C^*$-algebras.
This, actually, has been investigated by Blackadar and Kirchberg \cite {BK}
in their study of strong NF algebras.
We recall (by an equivalent definition from \cite {BK}) that
a $C^*$-algebra $\A$ is said to be a \emph{strong NF algebra}
if  for every  $x_1,\cdots, x_n \in \A$ and $\e > 0$,
there exists a  finite-dimensional $C^*$-algebra $\B$ and
a  (completely) isometric and complete order isomorphism
$\p$ from $\B$ onto a finite-dimensional $*$-subspace $F$ in
$\A$ such that
\[
dist(x_i, F)  < \e
\]
for all $i = 1, \cdots , n$.

It is clear from the definition that every strong NF algebra is a rigid
$\OL_{\infty, 1^+}$ space.
One of the major results in $\S 3$ is to show that these two notions are
actually  equivalent for unital $C^*$-algebras.
Indeed,  we prove in Theorem \ref {P2.inner} that if a unital (nuclear)
$C^*$-algebra is a rigid $\OL_{\infty, 1^+}$ space, then it is inner
quasi-diagonal, and thus  is a strong  NF algebra by \cite {BK2}.
We also prove in Theorem \ref {P2.quasi} that if a unital (nuclear)
$C^*$-algebra is an $\OL_{\infty, 1^+}$ space, then it is quasi-diagonal,
and thus is an NF algebra (see definition in \cite {BK}).

Summarizing our results in $\S 2$ and $\S 3$, we obtain the following
equivalent conditions for unital $C^*$-algebras.
\begin{theorem}
\label {P1.1}
Let $\A$ be a unital $C^*$-algebra. Then the following are equivalent:

\begin{itemize}
\item [(i)]  $\A$ is a strong  NF algebra (equivalently,
$\A$ is nuclear and inner quasi-diagonal),

\item [(ii)]  $\A$ is a  rigid ${\OL}_{\infty, 1^+}$ space,

\item [(iii)]  $\A$ is a  rigid rectangular ${\OL}_{\infty, 1^+}$ space,

\item [(iv)]  there exists a net of completely contractive projections
$P_\ga : \A \to \A$ such that $P_\ga \to id_\A$ in the point-norm topology.
\end{itemize}
\end{theorem}

Since every  rigid ${\OL}_{\infty, 1^+}$ space is
an ${\OL}_{\infty, 1^+}$ space, a $C^*$-algebra $\A$ is
an ${\OL}_{\infty, 1^+}$ space  if
it satisfies any of equivalent  conditions in Theorem \ref {P1.1}.
At  this moment,  we do not know whether ${\OL}_{\infty, 1^+}$  implies
rigid ${\OL}_{\infty, 1^+}$, and whether nuclearity and
quasi-diagonality  imply  ${\OL}_{\infty, 1^+}$ on (unital)
$C^*$-algebras.

During an operator space workshop organized by G. Pisier
at the IHP in Paris in January, 2000,
U. Haagerup showed the third author
that  if a unital $C^*$-algebra  satisfies the
condition (iv) in Theorem \ref {P1.1} then it is stably finite.
The strong NF algebra case was first proved by
Blackadar and Kirchberg  \cite {BK}.
Then it is natural to ask whether this is still true if $\A$ is an
$\OL_{\infty, 1^+}$ space, or an $\OL_{\infty, \la}$ space
for some $\la > 1$.
Along this line, we show in Theorem \ref {P.stablef} that
if a unital $C^*$-algebra  $\A$  is an $\OL_{\infty, \la}$ space
with $\la \le (\frac {1+{\sqrt 5}}2)^{\frac 12}$,
then $\A$ must be stably finite.

It is well-known from Banach space theory that if $\Omega$ is
a compact topological space, then the commutative $C^*$-algebra
$C(\Omega)$ is an $\SL_{\infty, 1^+}$ space and thus a rigid
$\SL_{\infty, 1^+}$ space  (see \cite {MP}).
Since it has  the \emph{minimal}  operator space structure,
it is also a rigid $\OL_{\infty, 1^+}$ space (see \cite {ER}),
and thus is a strong NF algebra by Theorem \ref {P1.1}.
This was also shown directly by Blackadar and Kirchberg \cite {BK}.
They actually  proved in \cite {BK} that
a quite large class of stably finite  nuclear $C^*$-algebras are
strong NF algebras and thus are (rigid) $\OL_{\infty, 1^+}$ spaces.
These  include  the spatial tensor product
$M_n \check \otimes C(\Omega)$,  their finite direct sums and
inductive limits  (such as AF algebras and AH algebras).
Moreover,  they  proved in \cite {BK2} that
subhomogeneous $C^*$-algebras,
and thus ASH algebras  are also  strong NF algebras.
A $C^*$-algebra is said to be {\em subhomogenuous}
if all of its irreducible representations are finite-dimensional
with  dim $\le n$ for some positive interger $n$,
and a $C^*$-algebra is said to be an {\em ASH}
algebra if it is the  inductive limit of subhomogenuous $C^*$-algebras.

In $\S 4$, we investigate the relation of the local structure of a nuclear
$C^*$-algebra and its second dual.
Using these results,  we are able to show in $\S 5$
that if $\A$ is a non-subhomogeneous nuclear $C^*$-algebra, then
$\A$ is an ${\OL}_{\infty, \lambda}$ space for every $\la > 6$
(see Theorem \ref {PNG.6}).
Combining the subhomogeneous and non-subhomogeneous cases, we obtain
the following theorem.

\begin{theorem}
\label {P1.2}
If $\A$ is a nuclear $C^*$-algebra, then $\A$ is an $\OL_{\infty, \la}$ space
for every $\la > 6$.
\end{theorem}

We  note that  Kirchberg proved in \cite {Ki2} that if $\A$ is a
separable non-type I nuclear $C^*$-algebra, then $\A$ is
completely isomorphic to the CAR algebra $\B$ with $d_{cb}(\A, \B) \le 256$,
and thus is an ${\OL}_{\infty, \la}$ space for every $\la > 256$
since  the CAR algebra is a (rigid) $\OL_{\infty, 1^+}$ space.
Our result significantly improves on the constant
obtained from Kirchberg's result.

 From these results, we see that in contrast to the $p=1$ case,
the $\OL_\infty$ constant $\la$  provides an interesting invariant
\[
{\OL}_{\infty} (\A) = \{ \la:  ~ \A ~{\rm is ~ an }~ \OL_{\infty, \la} ~
{\rm space} \}
\]
for nuclear $C^*$-algebras.
We can conclude  from Theorem \ref {P1.2} and Theorem \ref {P.stablef}
that if $\A$ is a nuclear $C^*$-algebra $\A$, then we have
\[
1 \le \OL_{\infty}(\A) \le 6,
\]
and if $\A$ is a nuclear non-stably finite unital $C^*$-algebra, we have
\[
(\frac {1+\sqrt 5} 2)^{\frac 1 2} < \OL_{\infty}(\A) \le 6.
\]
To end this paper,
we will make some remarks and propose some open  questions
related to this  new invariant $\OL_\infty$ on nuclear $C^*$-algebras
in $\S 6$.

Finally, we wish thank Bruce Blackadar, Uffe Haagerup,
Huaxin Lin, Eberhard Kirchberg
and Haskell Rosenthal for many stimulating discussions.

\section{Rigid rectangular ${\OL}_{\infty, 1^+}$ spaces}

Let us first consider the definition.
An operator space $V$ is said to be a {\em rigid rectangular}
${\OL}_{\infty, 1^+}$ \emph{space} if for every
$x _1, \cdots, x_n \in V$ and $\e > 0$,
there exits a finite direct sum of rectangular matrices
$\B =\oplus_{k=1}^l M_{m(k), n(k)}$  and a completely  isometric
injection $\p : \B \to V$ such that
\[
{\rm dist}(x_i, \p(\B)) < \e
\]
for all $i=1, \cdots, n$.

It is clear that $\B =\oplus_{k=1}^l M_{m(k), n(k)}$ can be identified
with the $(1, 2)$ off-diagonal corner of the finite-dimensional $C^*$-algebra
$\tilde \B = \oplus_{k=1}^l M_{m(k)+ n(k)}$, and thus is a finite-dimensional
ternary ring  of operators with the canonical ternary operation obtained
from $\tilde B$.
In general, an operator space $V \subseteq \B(K, H)$ is called
a {\em ternary ring of operators} (or simply, {\em TRO})
if it is   closed under the triple product
\[
(x, y, z) \in V \times V \times V \to xy^*z \in V.
\]
TORs were first introduced by Hestenes \cite {Hes} (see also \cite
{Harris}, \cite
{Zet}, \cite {KiTRO},  \cite {Hamana},  and \cite {EOR}).
A linear isomorphism between two TROs is a {\em TRO-isomorphism}
if it preserves the  triple products.
It is known that up to completely isometric TRO-isomorphism,
every  TRO can be identified  with the off-diagonal corner of a
unital $C^*$-algebra, and every finite-dimensional TRO has
the form $\oplus_{k=1}^l M_{m(k), n(k)}$ (see \cite {EOR}).

If $\A$ is a $C^*$-algebra and $P: \A \to \A$ is a
finite-rank completely contractive projection, then
it is known from Youngson \cite {Yo} that the range space
$P(\A)$ is a finite-dimensional TRO, and thus has the form
\[
P(\A) \cong \oplus_{k=1}^l M_{m(k), n(k)}.
\]

\begin{proposition}
\label {PR.1}
Let $\A$ be a $C^*$-algebra.
Then $\A$ is a rigid rectangular ${\OL}_{\infty, 1^+}$ space
if and only if there exists a net of finite-rank completely contractive
projections
$P_\ga : \A \to \A$ such that $P_\ga \to id_{\A}$ in the point-norm
topology.
\end{proposition}
\proof
Let us assume that  $\A$ is a rigid rectangular ${\OL}_{\infty, 1^+}$
space.
We let
\[
I = \{\ga = (x_1, \cdots, x_n, \e): ~ x_i \in \A, \e > 0 \}
\]
be the collection of finite subsets of $\A$ and $\e > 0$.
Given $\ga = (x_1, \cdots, x_n, \e) \in I$,
there exists a finite-dimensional TRO
$\B_\ga =\oplus_{k=1}^l M_{m(k), n(k)}$
and a completely  isometric injection
\[
\p_\ga : \B_\ga \to \A
\]
such that
\[
{\rm dist}(x_i, \p_\ga(\B_\ga)) < \e
\]
for all $i= 1, \cdots, n$.
Since  $\p_\ga(\B_\ga)$ is an injective subspace of $\A$,
the identity mapping on $\p_\ga(\B_\ga)$  has a completely contractive
extension $P_\ga: \A \to \p_\ga(\B_\ga)$,
which is a  projection from $\A$ onto $\p_\ga(\B_\ga)$.
In this case, there exist $w_1, \cdots, w_n \in \p_\ga(\B_\ga)$
such that $\|v_i-w_i\| < \e$ for all $i = 1, \cdots, n$.
Since $P_\ga(w_i) = w_i$, we have
\[
\|P_\ga(v_i) - v_i\| \le \|P_\ga (v_i - w_i)\| + \|w_i - v_i\| \le 2 \e.
\]
Then $\{P_\ga\}_{\ga \in I}$ (with the canonical partial order on
the  index set  $I$) is a net of  completely contractive projections on $\A$
such that  $\|P_\ga(x) - x\| \to 0$ for all $x \in \A$.

On the other hand, if we have a net of finite-rank completely
contractive
projections $P_\ga : \A \to \A$ such that
$\|P_\ga (x) - x\| \to 0$ for all $x \in \A$.
Then each  $P_\ga(\A)$ is a finite-dimensional TRO and thus is
completely isometric
to some $\oplus_{k=1}^l M_{m(k), n(k)}$.
Given   $x_1, \cdots, x_n \in \A$ and $\e > 0$,
there exists a completely contractive projection $P_\ga$ such that
$\|P_\ga(x_i) - x_i\| < \e$, and thus
\[
{\rm dist}(x_i, P_\ga(\A)) < \e
\]
for all $i = 1, \cdots, n$.
This shows that $\A$ is a rigid rectangular ${\OL}_{\infty, 1^+}$ space.
\endproof

The following proposition shows that if $\A$ is a unital $C^*$-algebra
and $P(1)$ is sufficiently close to the unital element $1$ of $\A$,
then $P(\A)$ must  be completely isometric to a finite-dimensional
$C^*$-algebra.

\begin{proposition}
\label {PR.3}
Let $\A$ be a unital $C^*$-algebra and let $P: \A \to \A$ be a finite-rank
completely contractive projection.
If $\|P(1) -  1\|_\A < {\frac 18}$, then
$P(\A)$ is completely isometric to a finite-dimensional $C^*$-algebra.
\end{proposition}
\proof
To simplify our notation, let us use  $V = P(\A)$ to denote the range
space of $P$.
Since $P$ is a finite-rank completely contractive projection on $\A$,
$V$ is a finite-dimesnional TRO with triple
product given by
\[
x\cdot y^*\cdot z = P(x y^* z)
\]
for all $x, y, z \in V$.
Up to (completely isometric) TRO-isomorphism,  we may write
\[
V = \oplus_{k=1}^lM_{m(k), n(k)},
\]
and identify $V$ with the off-diagonal corner (i.e. the $(1,2)$ corner)
of the finite-dimensional $C^*$-algebra $\B = \oplus_{k=1}^lM_{m(k)+n(k)}$.
Then
\[
C = \mbox{span}\{x \cdot y^*:  ~ ~ x, y \in V\}
\]
and
\[
D = \mbox{span}\{y^* \cdot z:  ~ ~  y, z \in V\}
\]
are  finite-dimensional $C^*$-subalgebras of $\B$,
and $V$ is a faithful non-degenerate $(C, D)$-bimodule.
The norms on $C$ and $D$ can be determined by the left module norm and
right module norm on $V$, respectively.
More precisely, for every $c \in C$ we have
\[
\|c\|_C = \mbox{sup}\{\|c \cdot x\|_V : ~  x\in V, \|x\|_V \le 1\},
\]
and for every $d \in D$ we have
\[
\|d\|_D = \mbox{sup}\{\|x \cdot d\|_V : ~  x\in V, \|x\|_V \le 1\}.
\]

If we let $a = P(1) \in V$, then $\|a\|_V = \|a\|_\A \le 1$, and
$a^*\cdot a$ is a positive element  in $D$ such that
\begin{eqnarray*}
\|a^*\cdot a\|_D & =&  \mbox{sup}\{\|x \cdot a^*\cdot a\|_V:
~ x \in V, ~\|x\|_V\le 1\}\\
&=& \mbox{sup}\{\|P(x a^*a)\|_V : ~ x \in V,~ \|x\|_V\le 1\}\\
&\le& \|a^*a\|_\A = \|a\|_\A ^2 \le 1.
\end{eqnarray*}
Moreover, we have
\begin{eqnarray*}
\|a^*\cdot a - 1_D\|_D &=& \mbox{sup}\{\|x\cdot a^* \cdot a - x\|_V:
x \in V, ~ \|x\|_V\le 1\}\\
&=& \mbox{sup}\{\|P(x  a^*  a - x)\|_V:
x \in V, ~ \|x\|_V\le 1\}\\
&\le & \mbox{sup}\{\|x(a^*  a - 1_\A)\|_\A:
x \in V, ~  \|x\|_V\le 1\}\\
&\le & \|a^*  a - 1\|_\A \le \|a^* - 1\|_\A \|a\|_\A +
\|a-1\|_\A < {\frac 14}.
\end{eqnarray*}
This shows that $a^*\cdot a$ is an invertible element in $D$, and
\[
{\frac 34} 1_D < a^*\cdot a \le 1_D.
\]
It follows that its square root $|a|$ in $D$ satisfies
\[
{\frac {\sqrt 3}  2} 1_D< |a| \le 1_D ,
\]
and thus
\[
\|1_D - |a|\, \|_D < 1 - {\frac {\sqrt 3}2 } < {\frac 14}.
\]
Regarding $a$ as an element in $\B$, we can consider the polar decomposition
$a = v \cdot |a|$ for some patial isometry $v \in \B$.
Since $|a|$ is an invertible element in $D$, and $V$ is a right
$D$-module,
we can conclude that
\[
v = a \cdot |a|^{-1} \in V,
\]
and thus  $v^*\cdot v$ is a projection in $D$.
We claim that $v^*\cdot v = 1_D$.

For any $x \in V$ with $\|x\|_V \le 1$,  we have
\begin{eqnarray*}
\|x \cdot v^*\cdot v - x \cdot a^* \cdot a \|_V &=& \|P(x v^*v -
xa^*a)\|_V
\le  \|v^*v - a^*a\|_\A\\
&\le & \|v^* - a^*\|_\A \|v\|_\A + \|a^*\|_\A \|v - a\|_\A  \\
&\le&   2 \|v - a\|_V \le 2 \|1_D - |a|\,\|_D < {\frac 12},
\end{eqnarray*}
and thus
\[
\|v^*\cdot v - a^* \cdot a \|_D < {\frac 12}.
\]
It follows that
\[
\|v^*\cdot v - 1_D\|_D \le \|v^*\cdot v - a^* \cdot a \|_D + \|a^* \cdot
a -1_D\|_D
\le {\frac 12} + {\frac 14} < 1.
\]
Since $v^*\cdot v$ is a projection in $D$,  we must have $v^*\cdot v =
1_D$.
Similarly, we can prove that $v\cdot v^* = 1_C$ in $C$.

In this case, there  is a natural $C^*$-algebra structure on $V$ given
by
\[
x\circ y = x \cdot v^* \cdot y ~ \, ~ \mbox{and} ~ \, ~ x^\sharp = v
\cdot x^* \cdot v
\]
for all $x, y \in V$ (see Zettl \cite {Zet}).
With this $C^*$-algebra structure, $v$ is the unit element of $V$.
The mapping
\[
\theta_D : x \in V \to v^* \cdot x \in D
\]
is a  unital $*$-isomorphism from $V$ onto  $D$.
The original matrix norm on $V$ coincides with the $C^*$-algebra
matrix norm on $V$ since $\theta_D$ is clearly a complete isometry
from $V$ onto $D$.
This completes the proof.
\endproof

It is clear that every rigid  ${\OL}_{\infty, 1^+}$ space is a rigid
rectangular
${\OL}_{\infty, 1^+}$ space.
But the converse is not necessarily true since for $n \ge 2$,
$M_{1, n}$ is clearly a rigid rectangular   ${\OL}_{\infty, 1^+}$ space,
but not a rigid ${\OL}_{\infty, 1^+}$ space (see \cite [$\S 4$]{ER}).
The following theorem shows that
the two notions coincide on unital $C^*$-algebras.

\begin{theorem}
\label{PR.2}
Let $\A$ be a unital $C^*$-algebra.  Then the following are equivalent:

\begin{itemize}
\item[{\rm (i)}]  $\A$ is a rigid  ${\OL}_{\infty, 1^+}$ space;

\item[{\rm (ii)}]  $\A$ is a rigid  rectangular ${\OL}_{\infty, 1^+}$
space;

\item[{\rm (iii)}]  there exists a net of finite-rank completely
contractive
projections $P_\ga : \A \to \A$ such that $P_\ga \to id_{\A}$ in the
point-norm topology.

\end{itemize}
\end{theorem}
\proof
It is obvious that  (i) $\Rightarrow$ (ii).
The equivalence  (ii)  $\Leftrightarrow$ (iii) is given by  Propositon
\ref {PR.1}.
If we have (iii), then there exists $\ga_0$ such that
\[
\|P_\ga(1) -1\| < {\frac 18}
\]
for all $\ga \ge \ga_0$.
It follows from  Proposition \ref {PR.3} that
$P_\ga(\A)$ are completely isometric to finite-dimensional
$C^*$-algebras for  all $\ga \ge \ga_0$.
Therefore, $\A$ is a rigid ${\OL}_{\infty, 1^+}$ space.
\endproof

\section{Quasi-diagonality  and inner quasi-diagonality}

Let $\A$ be a unital $C^*$-algebra and let $\psi: \A \to \B(H)$ be a unital
completely positive mapping.
It is known from the Stinespring
representation theorem that
there exists a Hilbert space $K$, a unital representation
(i.e. a unital $*$-homomorphism) $\pi: \A \to \B(K)$, and an isometry
$V \in \B(H, K)$ such that
\[
\psi (x) = V^*\pi (x) V ~\, ~\, ~ \, ~\, ~\, ~ \,~\, ~\, ~ \,
  ~\, ~\, ~ \, ~\, ~\, ~ \,~\, ~\, ~ \, (x \in \A).
\]
If we identify $H$ with the closed subspace $V(H)$ in $K$ and let $p$ be
the orthogonal projection from $K$ onto $H$, then we may regard $\psi$ as
the \emph{compression} of the unital representation $\pi$ by $p$ and write
\begin{equation}
\label {F2.comp}
\psi(x) = p \pi(x) p  ~\, ~\, ~ \, ~\, ~\, ~ \,~\, ~\, ~ \,
~\, ~\, ~ \, ~\, ~\, ~ \,~\, ~\, ~ \, (x \in \A).
\end{equation}
 From this, we may  easily  obtain the Schwarz inequality
\begin{equation}
\label {F2.Sch1}
\psi(x)^*\psi(x) \le \psi(x^* x)
\end{equation}
or equivalently,
\begin{equation}
\label {F2.Sch2}
0 \le \psi(x^* x)  - \psi(x)^*\psi(x)
\end{equation}
for all $x \in \A$.
Choi \cite {Choi} proved that for any $x\in \A$
if $\psi(x^*x) = \psi(x)^* \psi(x)$, then
\[
\psi(yx) = \psi(y)\psi(x)
\]
for all $y \in \A$.
In this case, we say that $x$ is in the {\em multiplicative domain} $D_\psi$ of
$\psi$.
The following is a very useful quantitative estimate of Choi's argument.
We thank the  referee for pointing out this simpler argument to us.

\begin{lemma}
\label{Pla1}
Let $\A$ and $\B$ be unital $C^*$-algebras and let $\psi: \A \to \B$
be a unital completely positive mapping. Then we have
\begin{equation}
\|\psi(x^*y) - \psi(x)^* \psi(y)\|
\le \|\psi(x ^* x) - \psi(x)^*\psi(x)\|^{\frac 12}
\|\psi(y^* y) - \psi(y)^*\psi(y)\|^{\frac 12}   \label {F2.1a}
\end{equation}
for all $x, y \in \A$.
\end{lemma}

\proof
We may first assume that $\B$ is a unital $C^*$-algebra on some Hilbert space,
and thus as we discussed in  (\ref {F2.comp}),
there exists a  Hilbert space $K$, a unital representation $\pi : \A \to \B(K)$
and an orthogonal projection $p$ on $K$ such that
\[
\psi(x) = p\pi(x)p
\]
for all $x \in \A$.
If we let $t: \A \to \B(K)$ be a complete contraction defined by
\[
t(x) = (1-p)\pi(x)p,
\]
then we have
\begin{equation}
\label {F3.t}
\psi(x^*y) - \psi(x)^* \psi(y) = p \pi(x^* y) p - p\pi(x^*)p \pi(y) p
= t(x)^* t(y)
\end{equation}
for all $x, y \in \A$.
It follows from (\ref {F3.t}) that
\begin{eqnarray*}
\|\psi(x^*y) - \psi(x)^* \psi(y)\|^2  &=&  \|t(x)^* t(y)\|^2
= \|t(y)^* t(x) t(x)^* t(y)\| \\
&\le& \|t(x) t(x)^*\| \|t(y)^* t(y)\|= \|t(x)^*t(x)\|\|t(y)^* t(y)\|\\
&=& \|\psi(x ^* x) - \psi(x)^*\psi(x)\|
\|\psi(y^* y) - \psi(y)^*\psi(y)\|.
\end{eqnarray*}
This proves the inequality (\ref {F2.1a}).
%We may obtain (\ref {F2.1b}) by considering
%$\tilde x = x^*$ and $\tilde y = y^*$ in (\ref {F2.1a}).
\endproof

It is well-known that quasi-diagonality
is a very important property for $C^*$-algebras.
We recall from Voiculescu \cite {Vo1} that a $C^*$-algebra $\A$ is
\emph{quasi-diagonal} if for every $x_1, \cdots, x_n \in \A$ and
$\e > 0$, there is a
representation $\rho$ of $\A$ on a Hilbert space $H$ and a
finite-rank projection $p\in \B(H)$ such that
\begin{equation}
\label {F2.quasi1}
\|[\rho(x_i), p]\| < \e  ~ \, {\rm and} ~\, \|x_i\| < \|p\rho(x_i)p\| +\e.
\end{equation}
For a unital $C^*$-algebra $\A$, this is equivalent to saying that for every
$x_1, \cdots,  x_n \in \A$ and $0 < \e < 1$,  there exists a unital
completely positive mapping
$\hat \psi: \A \to \B$ for some finite-dimensional $C^*$-algebra $\B$ such that
\begin{equation}
\label {F2.quasi2}
\|\hat \psi(x_i x_j) - \hat  \psi(x_i)\hat \psi(x_j)\| < \e
~\, {\rm  and } ~\,
\|x_j\| \le \|\hat \psi(x_j)\| +\e.
\end{equation}

\begin{theorem}
\label {P2.quasi}
If a unital {\rm (}nuclear{\rm )} $C^*$-algebra $\A$ is an
$\OL_{\infty, 1^+}$ space,
then $\A$ is quasi-diagonal.
\end{theorem}
\proof
It suffices to show that for every  contractive elements  $x_1,
\cdots,  x_n \in \A$ (with $x_1 = 1$)
and $0 < \e < 1$, we can find a unital complete positive mapping
$\hat \psi: \A \to \B$ satisfying  (\ref {F2.quasi2}).
We first let  $\de$ be a positive number with $0< \de < {\frac {\e^2} {16}}$.
Then there exists a finite-dimensional subspace $F$ in $\A$ such that
$F$ contains
all $x_j$ ($1\le j \le n$) and there exists a linear
isomorphism $\p$ from a finite-dimensional $C^*$-algebra  $\B$ onto
$F$ such that
$\|\p\|_{cb}  \|\p^{-1}\|_{cb} < 1+ \frac {\de^2}2$.
Without loss of generality, we may assume
that
\[
\|\p\|_{cb}< 1+ \frac {\de^2}2 ~\, {\rm and} ~\, \|\p^{-1}\|_{cb} \le 1.
\]
Since $\B$ is injective,  $\p^{-1}: F\to \B$ has a completely contractive
extension $\psi : \A \to \B$.

In the  following, we first show that we may suitably choose $\p$ and $\psi$
such that  $\psi(1)$ is a  positive element in $\B$.
Since  $1 = x_1 \in F$,  $b = \p^{-1}(1) = \psi(1)$ is a contractive element
in $\B$ such that $\p(b) = 1$ in $\A$.
If we let  $c = (1 - b^* b)^{\frac 12} \in \B$,
then
\[
\left \|\begin{bmatrix} b \\ c \end{bmatrix}\right \|^2 = \|b^* b +
c^* c\| = 1,
\]
and thus
\[
\left  \|\begin{bmatrix} 1 \\ \p(c) \end{bmatrix}\right \|^2
\le \|\p\|_{cb}^2 \left  \|\begin{bmatrix} b \\ c\end{bmatrix}\right \|^2
<  (1+ \frac {\de^2}2)^2.
\]
 From this, we can conclude that
\begin{equation}
\label {F2.pp1}
1 + \p(c)^* \p(c) < (1+ \frac {\de^2}2)^2 = 1 + \de^2 + \frac {\de ^4}4,
\end{equation}
and thus
\begin{equation}
\label {F2.pp2}
\|1 - b^* b\| = \|c\|^2 = \|\psi \circ \p(c)\|^2 \le
\| \p(c)\|^2 < \de^2 + \frac {\de^4} 4  < 2 \de^2 .
\end{equation}
This shows that  $b$ is an invertible element in $\B$.
If we let $b = v|b|$  be  the polar decomposition of $b$, then
$v$ must be a unitary in $\B$.
Moreover, since  $1-|b|$ and $1+|b|$ commute in $\B$,
  we have
\begin{equation}
\label{F2.pp3}
\|1 - |b|\| \le \|(1 - |b|)(1 +|b|)\| = \|1 - b^*b\| < 2 \de^2.
\end{equation}

Now we can modify $\psi$ and $\p$ by the unitary $v$.
We let $\tilde \p : \B \to F$ and $\tilde \psi: \A \to \B$
be mappings defined by
\[
\tilde \p(a) = \p(v a)  ~ \, {\rm and} ~\,  \tilde \psi (x) =  v^* \psi (x)
\]
for all $a \in \B$ and $x \in \A$.
Then it is easy to see that $\tilde \p$ is a linear isomorphism from
$\B$ onto $F$
with $\|\tilde\p\|_{cb} = \|\p\|_{cb} < 1+\e$ and $\tilde \psi$ is a completely
contractive extension of $\tilde \p^{-1}$ from $\A$ onto $\B$   such that
$\tilde \psi(1) = v^* b = |b|$ is  positive  in $\B$.

Next we construct a unital completely positive mapping $\hat \psi : \A \to \B$,
which is ``sufficiently close'' to $\tilde \psi $ and thus satisfies the
quasi-diagonal condition (\ref {F2.quasi2}).
For this purpose, let us assume that $\B$ is a unital $C^*$-subalgebra of
some matrix algebra  $\B({\Bbb C}^k)$.
It follows from Paulsen \cite {Paulsen} that there exists a unital
representation
$\rho : \A\to \B(H)$ for some Hilbert space $H$ and isometies
$V, W : {\Bbb C}^k \to H$ such that
\[
\tilde \psi (x) = V^* \rho (x) W
\]
for all $x \in \A$.
Since
\[
V^* W = \tilde \psi(1) = |b| =  \tilde \psi(1) ^* = W^* V,
\]
we obtain from (\ref {F2.pp3}) that
\begin{eqnarray*}
\|V-W\|^2 &=& \|(V - W)^*(V - W)\|= \|V^* V + W^* W - V^*W- W^* V\| \\
&=&  2\|1 - |b|\| < 4\de^2 .
\end{eqnarray*}
This shows that $\|V - W\| \le 2 \de $.

Since $\B$ is a unital injective $C^*$-subalgebra in $\B({\Bbb C}^k)$,
there exists a (completely positive) conditional expectation
$P: \B({\Bbb C}^k) \to \B$.  Then $\hat \psi : \A \to \B$ defined by
\[
\hat \psi(x)  = P(V^* \rho (x) V)
\]
is a unital completely positive mapping  from $\A$ into $\B$.
Since the range of $\tilde \psi$ is contained in $\B$, we have
\[
\hat \psi - \tilde \psi = P(V^* \rho  V) - V^* \rho W = P(V^* \rho (V-W)),
\]
and thus
\begin{equation}
\label{F2.pp4}
\|\hat \psi - \tilde \psi\|_{cb} \le \|V -W\| < 2\de.
\end{equation}

Given  a unitary element $u \in \B$, if we let $y = \tilde \p (u)$,
then  $\|y\| < 1+\frac {\de^2}2$
and $\tilde \psi(y)  = \tilde \psi\circ \tilde \p(u)  = u$.
It follows that
\[
\|\hat \psi(y) - u\| =  \|\hat \psi(y) - \tilde \psi(y)\| < 2\de\|y\| < 3\de,
\]
and thus
\[
\|\hat \psi(y)^* \hat  \psi(y) - 1 \|
\le \|\hat \psi(y)^* \| \|\hat \psi(y)  - u\|
+\|\hat \psi(y)^* - u^* \| \|u\| < (2 + {\frac {\de^2}2})3\de < 8 \de.
\]
 From the  Schwarz inequality (\ref {F2.Sch1}), we get
\[
1 - 8 \de \le  \hat \psi(y)^* \hat \psi(y) \le
\hat \psi(y^* y) < (1+ {\frac{\de^2}2})^2,
\]
and thus
\begin{equation}
\label {F3.SS1}
0\le {\hat \psi}(y^* y)-  {\hat \psi}(y)^* {\hat \psi}(y)
\le 8  \de +  \de^2 +  {\frac{\de^4}4} < 16 \de.
\end{equation}
On the other hand, for every  $x \in \A$  we get the inequality
\begin{equation}
\label {F3.SS2}
0\le \hat \psi (x^* x) -  \hat \psi(x)^* \hat \psi(x) \le \hat \psi (x^* x)
\end{equation}
from  (\ref {F2.Sch2}).
Then applying  Lemma \ref {Pla1} together with (\ref  {F3.SS1}) and
(\ref  {F3.SS2}),  we obtain
\begin{eqnarray*}
\|\hat \psi(x y) - \hat \psi(x) \hat \psi(y)\|
&\le & \|\hat \psi(x ^* x) - \hat \psi(x)^*\hat \psi(x)\|^{\frac 12}
\|\hat \psi(y^* y) - \hat \psi(y)^*\hat \psi(y)\|^{\frac 12}  \\
&<& 4 {\sqrt \de} \|x\| < \e \|x\|
\end{eqnarray*}
for all $x \in \A$.

In general,  if we are given any  $y \in F$ with $\|y\| < 1$,
then $a = \tilde \psi (y)$ is contained in the open unit ball of $\B$, and
thus  can be written as a convex combination of unitary elements in
$\B$ by the Russo-Dye theorem, i.e. there exist unitary elements
$u_i \in \B$ and positive numbers $\al_i$ with $\sum _i \al_i = 1$
such that  $a =\sum \al_i u_i$.
In this case, we can write
\[
\tilde \p(a) = \sum \alpha_i \tilde \p(u_i) = \sum \alpha_i y_i
\]
as a convex combination of $y_i = \tilde \p(u_i)$, and thus obtain
\[
\|\hat \psi(xy) -\hat  \psi(x)\hat  \psi(y)\| \le
\sum_i \al_i \|\hat  \psi(xy_i) - \hat  \psi(x)\hat \psi(y_i)\|<   \e \|x\|.
\]
By the continuity of $\hat \psi$, we can conlcude that
\[
\|\hat \psi(xy) -\hat  \psi(x)\hat  \psi(y)\| \le
\e \|x\| \|y\|
\]
for all $x\in \A$ and $y \in F$.
This proves the first inequality in (\ref {F2.quasi2}).

The second inequality follows from
\begin{eqnarray*}
\|x_j\| &=& \|\tilde \p \circ \tilde \psi (x_j) \|
< (1+{\frac {\de^2}2}) \|\tilde \psi (x_j) \|  \\
&\le & (1+{\frac {\de^2}2})  (\|\hat \psi(x_j)\| + \|\hat  \psi (x_j) -
\tilde \psi (x_j) \|) \\
&<&  \|\hat \psi(x_j)\| + {\frac {\de^2}2} + 2 \de < \|\hat \psi(x_j)\| + \e.
\end{eqnarray*}
\endproof

Blackadar and Kirchberg studied NF algebras and strong NF algebras in
\cite {BK}.
They proved that  a unital $C^*$-algebra is an NF algebra (see definition given
in \cite {BK}) if and only if it is nuclear and quasi-diagonal.
Moreover, they characterized strong NF algebras by nuclearity and inner
quasi-diagonality  in \cite {BK2}, where the \emph{inner} condition requires
that  the finite-rank projection $p$ (in the definition of quaisi-diagonality)
is contained in $\pi(\A)''$.
More precisely,  a $C^*$-algebra $\A$ is said to be \emph{inner quasi-diagonal}
if  for every $x_1, \cdots , x_n \in \A$ and
$\e > 0$, there is a representation $\rho$ of $\A$ on a Hilbert space $K$
and a finite-rank projection $p \in \rho(\A)'' \subseteq \B(K)$ such that
(\ref {F2.quasi1}) is satisfied.
This is a stronger condition than quasi-diagonality
(see examples given in \cite {BK2}).

\begin{theorem}
\label {P2.inner}
If a unital {\rm (}nuclear{\rm )} $C^*$-algebra $\A$ is a rigid
$\OL_{\infty, 1^+}$ space,  then it is  inner quasi-diagonal,
and thus is a strong NF algebra.
\end{theorem}
\proof
Given  contractive elements
$x_1, \cdots, x_n \in \A$ (with $x_1 = 1$) and $0 < \e < 1$,
we let $\de$ be a positive number  with $0 < \de < \frac {\e^2}{16}$.
Since $\A$ is a rigid $\OL_{\infty, 1^+}$ space, there exists a complete
isometry $\p: \B \to F$ from a finite-dimensional $C^*$-algebra
$\B$ onto a finite-dimensional subspace $F$ of $\A$ such that
\[
dist(x_i, F) < {\frac{\de} 4}
\]
for all $1\le i\le n$.
Then for each $i$, we may find an element $a_i \in \B$ such that
\begin{equation}
\label{F3.ai}
\|\p(a_i) -x_i\| < {\frac{\de} 4}  ~\mbox{and} ~ \|\p(a_i)\| < 1+\frac \de 4.
\end{equation}
We note that in contrast to the argument in Theorem \ref {P2.quasi},
$x_1= 1$ need not be in $F$.   We can only approximate it by $\p(a_1)$.
Therefore, we need a modified argument given as follows.

We let $\psi: \A \to \B$ be  a completely contractive  extension of
$\p^{-1}$, and  let $b = \psi(1)$ and $c = (1 - b^* b)^{\frac 12} \in \B$.
Then we have
\[
\|a_1-b\|= \|\psi\circ \p(a_1) - \psi(1)\| \le \|\p(a_1) - 1\| < {\frac{\de} 4}
\]
and
\[
\left \|\begin{bmatrix} b \\ c \end{bmatrix}\right \| = \|b^* b + c^*
c\|^{\frac 12} = 1.
\]
It follows that
\[
\left \|\begin{bmatrix} 1 \\ \p(c) \end{bmatrix}\right \|^2
< \left (\left \|\begin{bmatrix} \p(a_1) \\ \p(c) \end{bmatrix}\right \|
+ {\frac{\de} 4} \right )^2
= \left ( \left \|\begin{bmatrix} a_1 \\ c\end{bmatrix}\right \|
+ {\frac{\de} 4} \right ) ^2
<  \left ( \left \|\begin{bmatrix} b \\ c\end{bmatrix}\right \| +
{\frac{\de} 2}\right ) ^2 = (1+ {\frac{\de} 2}) ^2.
\]
Using the same argument given for (\ref {F2.pp1}), (\ref {F2.pp2}) and
(\ref {F2.pp3}), we can prove that
$b$ is an invertible element in $\B$ and if $b = v |b|$  is the polar
decomposition for $b$, then   $v$ is a unitary in $\B$ and
\[
0 \le \|1  - |b|\| \le 2 \de^2.
\]
Then we may define a complete isometry $\tilde \p : \B \to F$ and
a complete contraction $\tilde \psi: \A \to \B$
by letting
\[
\tilde \p(a) = \p(v a) ~ \, {\rm and} ~\,  \tilde \psi (x) =  v^* \psi (x)
\]
for all $a \in \B$ and $x \in \A$.
It is easy to see that $\tilde \psi(1) = |b|$ is a positive element in
$\B$ and we have $\tilde \psi \circ \tilde \p = id_\B$ and
$\tilde \p \circ \tilde \psi_{|F} = id_F$.

We let
\[
\L_\B = \left \{\begin{bmatrix} \la & a  \\ b^* & \mu \end{bmatrix}:
~ \la, \mu \in {\Bbb C}, a, b \in \B \right \}
\]
denote  the operator system in $M_2(\B)$ induced by $\B$, and let
$\Phi: \L_\B \to M_2(\A)$ be the  unital completely
positive mapping  from  $\L_\B$ onto the operator system $\L_F$ in $M_2(\A)$,
which is defined by
\[
\Phi\left (\begin{bmatrix} \la & a  \\ b^* & \mu \end{bmatrix}\right )
= \begin{bmatrix} \la & \tilde \p(a)  \\ \tilde \p(b)^* &\mu\end{bmatrix}.
\]
We let
\[
\L_\A = \left \{\begin{bmatrix} \la & x  \\ y^* & \mu\end{bmatrix} :
\la, \mu \in
{\Bbb C}, x, y \in \A \right \}
\]
denote the operator system in $M_2(\A)$ induced by $\A$, and let
$\Psi: \L_\A \to M_2(\B)$ be the unital completely positive mapping defined by
\[
\Psi\left ( \begin{bmatrix} \la & x  \\ y^* & \mu \end{bmatrix}\right )
= \begin{bmatrix} \la & \tilde \psi(x)  \\ \tilde \psi(y)^* &\mu\end{bmatrix}.
\]
Then $\Psi$ extends to a unital completely positive mapping, which is
still denoted
by $\Psi$, from $M_2(\A)$ into $M_2(\B)$.  It is easy to see that
$\Psi \circ \Phi = id_{\L_\B}$, and thus $\Phi$ is a completely order
isomorphism
from $\L_\B$ onto $\L_F$.
We note that  $\L_F$ is contained in the multiplicative domain $D_\Psi$ of
$\Psi$, i.e. for every  $\tilde y \in \L_F$, we have
\[
\Psi(\tilde x \tilde y) = \Psi(\tilde x) \Psi(\tilde y)
\]
for all $\tilde x \in M_2(A)$.
To see this, given any $\tilde a = \begin{bmatrix} \la & u  \\ v^* & \mu
\end{bmatrix}  \in \L_\B$ with $u, v$ being unitaries in $\B$,
$\tilde a^* \tilde a$ is again an element in $\L_\B$ and thus satisfies
\[
\tilde a^* \tilde a = \Psi (\Phi(\tilde a)^*)\Psi(\Phi(\tilde a))
\le \Psi (\Phi(\tilde a^*)\Phi(\tilde a)) \le \Psi (\Phi(\tilde a^*\tilde a))
= \tilde a^* \tilde a.
\]
Therefore, we have
\[
\Psi (\Phi(\tilde a)^*)\Psi(\Phi(\tilde a)) = \Psi (\Phi(\tilde
a)^*\Phi(\tilde a))
\]
and thus $\tilde y =\Phi(\tilde a)$ is contained in the  multiplicative domain
of $\Psi$  by Choi \cite {Choi}.
Since the open unit ball of  $\B$ is contained in convex hull of
unitary elements in
$\B$, we may conclude that $\L_F = \Phi(\L_\B)$ is contained in the
multiplicative
domain $D_\Psi$ of $\Psi$.

We let $\A_\Psi$ denote the unital subalgebra of $\A$ generated by $\L_F$.
It is easy to see that $\A_\Psi$ is contained in the  multiplicative
domain  $D_\Psi$.
Since $\L_F$ is an operator system in $\A$, $\A_\Psi $ must be
self-adjoint and thus is a unital $C^*$-subalgebra of $\A$.
Then $\Psi$ restricted to $\A_\Psi$ induces a unital
$*$-homomorphism
$\pi = \Psi_{|\A_\Psi}$ from $\A_\Psi$ into  $M_2(\B)$.
Since $\L_\B$ is contained in the range of $\pi$, it is easy to see that $\pi$
maps $\A_\Psi$ onto $M_2(\B)$.
We note that this is a generalization of an argument in Choi and
Effros \cite [Theorem 4.1]{CE0}.

We may write $\B = \oplus_k \B({\Bbb C}^{n_k})$ and $\pi = \oplus_k \pi_k$,
where each  $\pi_k$ is an irreducible representation from $\A_\Psi$ onto
$M_2(\B({\Bbb C}^{n_k})) = \B({\Bbb C}^{2n_k})$.
For each $k$, we may extend $\pi_k$ to an irreducible representation
$\tilde \pi_k : M_2(\A) \to \B(\tilde H_k)$ on a larger Hilbert space
$\tilde H_k$.
Up to unitary equivalence, we may write
\[
\tilde \pi_k = id_2 \otimes \rho_k,
\]
where $\rho_k: \A \to \B(K_k)$ is an irreducible representation,
and we can write
\[
\pi_k = \begin{bmatrix} V_k^* & 0 \\ 0 &W_k^* \end{bmatrix}
(id_2 \otimes \rho_k)
\begin{bmatrix} V_k & 0 \\ 0 &W_k \end{bmatrix}
\]
for some isometries $V_k$ and $W_k$ from
${\Bbb C}^{n_k}$ into $K_k$.
If we let $V = \oplus_{k=1}^n V_k$ and let $W = \oplus_{k=1}^n  W_k$, then
$V$ and $W$ are isometries satisfying
\[
V^* W = V^* \rho(1) W = \tilde \psi (1) \ge 0
\]
in $\B$.
It follows (by the same argument given in Theorem \ref {P2.quasi})
that we have
\[
\|V_k - W_k\| \le \|V-W\| < 2\de,
\]
and the orthogonal projections $p_k = V_k V^*_k$ and $q_k = W_k W^*_k$
from  $K_k$ on ${\Bbb C}^{n_k}$ satisfy
\begin{equation}
\label {F3.pqk}
\|p_k - q_k\| \le \|V_k - W_k\|\|V^*_k\| + \|W_k \|V^*_k - W^*_k\| < 4\de
\end{equation}
for all $k$.
Since the range of $\tilde \pi_{k |\A_\Psi} = \pi_k$ is contained
in $M_2(\B({\Bbb C}^{n_k})) = \B({\Bbb C}^{n_k}\oplus {\Bbb C}^{n_k})$,
the orthogonal  projection  $p_k \oplus q_k =
\begin{bmatrix} p_k & 0 \\ 0 & q_k
\end{bmatrix}$   from $K_k \oplus K_k$ onto  ${\Bbb C}^{n_k}\oplus
{\Bbb C}^{n_k}$
leaves $\tilde \pi_k(\A_\Psi)$ invariant, i.e. we have
\[
\tilde \pi_k (\tilde x) (p_k\oplus q_k)  = (p_k\oplus q_k)  \tilde
\pi_k (\tilde x)
\]
for all $\tilde x \in \A_\Psi$.
Then for any $y \in F$,
\[
\tilde \pi_k\left ( \begin{bmatrix} 0 & y \\ 0 & 0 \end{bmatrix}\right )
\begin{bmatrix} p_k & 0\\ 0 & q_k \end{bmatrix}
= \begin{bmatrix} p_k & 0\\ 0 & q_k \end{bmatrix}
\tilde \pi_k\left ( \begin{bmatrix} 0
& y \\ 0 & 0 \end{bmatrix}\right )
\]
implies
\begin{equation}
\label {F3.pqky}
p_k \rho_j(y) = \rho_j(y) q_k.
\end{equation}

In general, $\rho_k$ are not necessarily all non-equivalent.
We let $\rho_1, \cdots, \rho_r$ be non-equivalent irreducible representations,
and for $k> r$, we may choose $l\le r$ such that $\rho _k \cong \rho_l$ and
identify $K_k$ with $K_l$, and $V_k$, $W_k$ with the corresponding
isometries $V_l, W_l$ from ${\Bbb C}^{n_l}$ into  $K_l$.
If we let  $\rho = \oplus_{k=1}^r \rho_k$,  $K = \oplus_{k=1}^r K_k$,
then
\[
p = \oplus_{k=1}^r p_k ~\, {\rm and } ~\, q = \oplus_{k=1}^r q_k
\]
are  projections  contained in
$\oplus_{k=1}^r \B({\Bbb C}^{n_k}) = \rho(\A)''$, and
we can conclude from (\ref {F3.pqk}) and (\ref {F3.pqky}) that
\begin{equation}
\label{F3.pq}
\|p - q\| = \mbox{max} \{\|p_k-q_k\|\} < 4\de
\end{equation}
and
\begin{equation}
\label{F3.pqy}
p\rho(y) = \rho(y) q
\end{equation}
for all $y \in F$.
If we let  $y_i = \p(a_i) \in F$, then we have from (\ref {F3.ai}),
(\ref {F3.pq}) and (\ref {F3.pqy}) that
\begin{eqnarray*}
\|\rho (x_i)p - p\rho(x_i)\| &\le& \|\rho (y_i)p - p\rho(y_i)\| + 2
\|x_i-y_i\|\\
&\le &  \|p-q\|\|y_i\| + 2 \|x_i-y_i\|\\
& < & 4\de (1+ {\frac {\de} 4})  +  {\frac{\de} 2}  < \e,
\end{eqnarray*}
and
\begin{eqnarray*}
\|x_i\| &<& \|y_i\| + \frac {\de} 4 = \|\tilde \psi(y_i)\|  +\frac {\de} 4  \\
&\le& \|p\rho(y_i)q\| + \frac {\de} 4 <
\|p\rho(y_i)p\| + \|p-q\|\|y_i\| +\frac {\de} 4 \\
&<& \|p\rho(x_i)p\| + \e.
\end{eqnarray*}
This shows that $\A$ is inner quasi-diagonal.
\endproof

Summarizing our results in $\S 2$ and $\S 3$, we have obtained the following
relations for unital $C^*$-algebras:
\[
\mbox{strong ~ NF } \Leftrightarrow \mbox{rigid} {\OL}_{\infty, 1^+}
\Leftrightarrow \mbox{rigid rectangular } {\OL}_{\infty, 1^+}\Rightarrow
{\OL}_{\infty, 1^+} \Rightarrow {NF}.
\]

Recall that a unital $C^*$-algebra is said to be {\em finite }
if every isometry must be a unitary.
This is equivalent to saying that if $v $ is a partial isometry in $\A$,
then $p = v^* v \le q = v v^*$ implies that $p = q$.
A unital $C^*$-algebra $\A$ is said to be {\em stably finite} if for
every
$n \in \N$,  $M_n(\A)$ is finite.

\begin{theorem}
\label{P.stablef}
Let $\A$ be a unital $C^*$-algebra.
If  $\A$ is an $\mathcal{OL}_{\infty,\lambda}$ space
with $\lambda\le (\frac{1+ \sqrt 5}2)^{\frac 12}$,
then $\A$ must be a stably finite $C^*$-algebra.
\end{theorem}

\begin{proof}
We first note that the $\OL_{\infty, \lambda}$ structure is a stable property
on operator spaces, i.e. if $\A$ is an $\mathcal{OL}_{\infty,\lambda}$ space,
then   for each $n \in {\Bbb N}$,   $M_n(\A)$ is again an $\OL_{\infty, \la}$
space.
Therefore, it surffices to show that $\A$ is a finite $C^*$-algebra.
The stable finiteness follows immediately.

Let us assume that $\A$ is not finite. Then there is an isometry
$s \in \A$ such that  $p=1-s s^*\neq 0$.
Since $\A$ is an $\mathcal{OL}_{\infty,\lambda}$ space,
there exists a finite-dimensional subspace $F \subseteq \A$ which
contains $\{s, p\}$ and there exists a linear isomorphism
\[
\p: F \to \B
\]
from $F$ onto a finite-dimensional $C^*$-algebra $\B$ such that
\[
\|\p\|_{cb}< 1  ~\, {\rm and} ~ \, \|\p^{-1}\|_{cb}\le \lambda.
\]
Since  $\|[s\, ~ p]\| = \|s s^* + (1 - s s^*)\| = 1$ and
${\frac 1 \lambda} = {\frac 1 \lambda} \|p\| < \|\p (p)\|$,
we  have
\[
\|[\p(s)\ \p (p)]\|\le\|[s\ p]\|=1
\]
and thus
\[
{\frac 1 {\lambda^2}} <  \|\p (p)\|^2 =
\|\p(p)\p(p)^*\| \le \|1-\p(s)\p(s)^*\|.
\]
On the other hand,  let us consider $b=(1-\p(s)^* \p(s))^{\frac 12}\in \B$.
Since $s^*s=1$ and
\[
\|1+ \p^{-1}(b)^* \p^{-1}(b)\|^{\frac 12} =
\left \|\begin{bmatrix} s \\ \p^{-1}(b) \end{bmatrix}\right \|
\le\|\p^{-1}\|_{cb}\left \|\begin{bmatrix} \p(s) \\ b \end{bmatrix}\right \|
<\lambda,
\]
we have
\[
\|1-\p(s)^*\p(s)\| = \|b\|^2 \le \|\p^{-1}(b)\|^2 = \|\p^{-1}(b)^*\p^{-1}(b)\|
<\lambda^2 -1
\le \frac {\sqrt 5 - 1} 2 < 1.
\]
Then $\p(s)^*\p(s)$ and thus $\p(s)$ must be invertible elements in
the finite-dimensional $C^*$-algebra $\B$.
Considering the polar decomposition of $\p(s)$, we get
\[
\|1-\p(s)\p(s)^*\| = \|1-\p(s)^*\p(s)\|,
\]
and thus
\[
\frac 1 {\lambda^{2}}\le\|1-\p(s)\p(s)^*\|
=\|1-\p(s)^*\p(s)\|<\lambda^2-1 .
\]
This implies $\lambda > (\frac {1+ \sqrt 5} 2)^{\frac 1 2}$
and leads to a contradiction  to the hypothesis.
\end{proof}

\section{${\OL}_{\infty, \la}$ structure  related to the second duals}

Let us first recall that an operator space $V$ is said to be
{\em locally reflexive} if for any finite-dimensional operator space
$\tilde F,$
each complete contraction $\varphi: \tilde F\rightarrow V^{**}$ can be
approximated in the point-weak$^{*}$ topology by complete contractions
$\psi : \tilde F\rightarrow V$.
Equivalently, given any finite-dimensional subspace  $G\subseteq V^{*}$
and $\varepsilon >0,$
we can find a mapping $\psi : \tilde F\rightarrow V$ such
that
\begin{equation}  \label{FG.1}
\left\| \psi \right\| _{cb}<1+\varepsilon \text{ and }\langle
\psi(x),f\rangle =\langle \varphi (x),f\rangle
\end{equation}
for all $x
\in \tilde F$ and $f\in G$ (see \cite [Lemma 6.4]{EJR}).
An operator space $V$ is said to be \emph{strongly locally reflexive} if
given any finite-dimensional subspaces $\tilde F\subseteq V^{**}$,
$G\subseteq V^{*}$
and $\e >0$, there exists a linear isomorphism $\psi$ from $\tilde F$ onto
a subspace
$F$ of $V$ such that
\begin{enumerate}
\item[(a)]  $\left\| \psi \right\| _{cb}\left\| \psi ^{-1}\right\|
_{cb}<1+\epsilon,$

\item[(b)]  $\langle \psi (v),\, f\rangle =\langle v,\,  f\rangle $ for
all
$v\in \tilde F$ and $f\in G,$ and

\item[(c)]  $\psi (v)=v$ for all $v\in \tilde F\cap V.$
\end{enumerate}
It is obvious that every strongly locally reflexive operator space is
locally reflexive, and it was shown in \cite {EJR} that the operator
preduals
of von Neumann algebras are {\em all} strongly locally reflexive.

\begin{lemma}
\label {PG.1}
Suppose that $V$ and $W$ are finite-dimensional operator spaces
and ${\A}$ is a $C^*$-algebra.
If a linear mapping $\varphi: V \to W$ has a completely bounded
factorization
\[
\begin{array}{ccccc}
&  & {\A}^{**} &  &  \\
& {\scriptstyle \al }\longnearrow  &  & \longsearrow {\scriptstyle
\be } &  \\
{V} \hspace{-10 pt} &  & \stackrel{\varphi}{\hspace{-10pt}{\hbox to 30pt
{\rightarrowfill}}\hspace{-10pt} } &  &\hspace{-10 pt} W
\end{array} ,
\]
then for any $\e > 0$   we may replace $\be$ in above
diagram by a weak$^*$ continuous
completely bounded mapping $\tilde \be : {\A}^{**} \to W$ such that
$\|\tilde \be\|_{cb} \le (1+\e) \|\be\|_{cb}$.
\end{lemma}
\proof
Let $\be^* : W^* \to {\A}^{***}$ be the adjoint mapping of $\be$.
Then $\be^*(W^*)$  is a finite-dimensional subspace of ${\A}^{***}$.
Since  the operator dual ${\A}^*$ of the
$C^*$-algebra ${\A}$ is  locally reflexive (see \cite{EJR}),
for any $\e > 0$ and the finite-dimensional subspace $\al(V) \subseteq
{\A}^{**}$, there exists a  mapping  $\psi :  \be^*(W^*) \to {\A}^*$
such that   $\|\psi\|_{cb} < 1+\e$ and
\[
\langle \al(v),\, \psi \circ \be^*(f) \rangle  =  \langle \al(v) , \,
\be^*(f)\rangle
\]
for all $f \in W^*$ and $v \in V$.
Then the adjoint mapping
\[
\tilde \be = (\psi \circ \be^* )^*  : {\A}^{**} \to W
\]
is  weak$^*$ continuous  such that
\[
\|\tilde \be\|_{cb} \le \|\psi\|_{cb} \|\be^*\|_{cb} <
(1+\e)\|\be\|_{cb},
\]
and
\[
\langle \tilde \be \circ \al (v), f\rangle
= \langle \al(v), \psi \circ \be^*(f) \rangle
=  \langle \al(v), \be^*(f) \rangle =  \langle \be\circ \al(v), f
\rangle
= \langle \varphi(v), f \rangle
\]
for all $f \in W^*$ and $v \in V$. From this we can conclude that
$\varphi = \tilde \be \circ \al$.
\endproof

The following lemma is a generalization of \cite [Theorem 7.6]{EJR}
from the case $\A = K(H)$  to general locally reflexive $C^*$-algebras.

\begin{lemma}
\label{PG.2}
Suppose that  $V$ and $W$ are  finite-dimensional operator spaces and
$\A$ is a locally reflexive $C^*$-algebra.
If  a linear mapping $\varphi : V \to W$ has
a completely bounded factorization
\begin{equation}
\label {FG.2}
\begin{array}{ccccc}
&  & {\A}^{**} &  &  \\
& {\scriptstyle \al }\longnearrow  &  & \longsearrow {\scriptstyle
\be } &  \\
{V} \hspace{-10 pt} &  & \stackrel{\varphi}{\hspace{-10pt}{\hbox to 30pt
{\rightarrowfill}}\hspace{-10pt} } &  &\hspace{-10 pt} W
\end{array} ,
\end{equation}
then for any $0< \e < 1$ we may obtain a completely bounded
factorization
\begin{equation}
\label {FG.3}
\begin{array}{ccccc}
&  & {\A} &  &  \\
& {\scriptstyle \tilde \al }\longnearrow  &  & \longsearrow
{\scriptstyle
\tilde \be } &  \\
{V} \hspace{-10 pt} &  & \stackrel{\varphi}{\hspace{-10pt}{\hbox to 30pt
{\rightarrowfill}}\hspace{-10pt} } &  &\hspace{-10 pt} W
\end{array}
\end{equation}
such that
\begin{equation}
\label {F.addition}
\|\tilde \al\|_{cb} < (1+\e)\|\al\|_{cb} ~\mbox{and}~
\|\tilde \be\|_{cb}< (1+\e)\|\be\|_{cb}.
\end{equation}
Moreover, we can choose $\tilde \alpha$ such that
\begin{equation}
\label {F.strongloc}
\al(V)\cap \A \subseteq \tilde \al(V).
\end{equation}
\end{lemma}
\proof
Given $0 < \e < 1$, it follows from Lemma \ref {PG.1} that
we may first replace
$\be$ in (\ref {FG.2}) by a weak$^*$ continuous mapping
$\tilde \be: {\A}^{**} \to W$
such that
\[
  \varphi = \tilde \be \circ \al ~ \mbox{ and} ~
\|\tilde \be\|_{cb} < (1+\e)\|\be\|_{cb}.
\]
The pre-adjoint  $\tilde \be_{*}: W^{*}\rightarrow \A^{*}$ is a
well-defined mapping from $W^*$ into $\A^*$, and
$\tilde \be_*(W^*)$ is a finite-dimensional subspace of $\A^*$.
We let $\tilde F =\al(V)$ denote the finite-dimensional subspace
in  $\A^{**}$, and let $n = {\rm dim} \tilde F $ and $0<\delta \le
\frac\e {2 n^{3}}$.
Since $\A$ is locally reflexive,
there exists a mapping $\psi: \tilde F \rightarrow \A$ such
that
\begin{equation}
\left\| \psi \right\| _{cb}<1+\delta  < 1+\varepsilon,  \label{E3a}
\end{equation}
and
\[
\langle \psi \circ \al (v), \tilde \be_*(f)\rangle
=\langle \al (v), \tilde \be _*(f)\rangle
\]
for all $v\in V$ and $f\in W^{*}$.
This shows that
\begin{equation}
\tilde \be \circ \psi \circ \al  = \tilde \be \circ \al = \varphi,
\label {E4}
\end{equation}
and we obtain (\ref {FG.3}) and (\ref {F.addition})
by letting $\tilde \alpha = \psi\circ \alpha$.

Note that in general an infinite-dimensional $C^*$-algebra
(for example,  $K(\ell_2)$) is not strongly locally
reflexive.
So we need  the following argument to obtain $\tilde \alpha$ satisfying
(\ref{F.strongloc}).
We let $C\subseteq \mathcal{CB}(\tilde F,\A)$ be the
convex set of all mappings $\psi: \tilde F\rightarrow \A$
satisfying (\ref{E3a}) and
(\ref{E4}). We let $F= \tilde F\cap \A$ and let $\iota: F\rightarrow \A$ be
the inclusion mapping.
We let $C_{0}\subseteq \mathcal{CB}(F,\A)$
denote the convex set of all mappings $\psi \circ\iota,$
where $\psi \in C$.
We claim that $\iota$ is in the point-norm closure of $C_{0}.$
This is apparent since if we are given an arbitrary finite-dimensional
subspace $G\subseteq \A^{*},$ then our previous argument shows that
there is a mapping $\psi ^{\prime }: \tilde F\rightarrow \A$
satisfying
\[
\left\| \psi ^{\prime }\right\| _{cb}<1+\delta ,
\]
and
\[
\langle \psi ^{\prime }\circ \al (v), f\rangle =\langle \al (v),
f\rangle
\]
for all $v\in V$ and $f\in \tilde \be_{*}(W^{*})\,+G$.
The latter equation implies that
\[
\tilde \be \circ \psi^\prime \circ \al = \tilde \be \circ  \al =
\varphi.
\]
This shows that  we can suitably choose
a net of $\psi ^{\prime }$ in $C$ such that $\psi^{\prime }\circ \iota \in C_0$
converges to $\iota$ in the point-weak topology.
Then the usual convexity argument shows that $\iota$
is in the point-norm closure of $C_0$.
Since $F$ is finite-dimensional, its closed unit
ball is totally bounded and thus
we may choose a mapping $\psi \in C\ $ such that
\[
\left\| \iota -\psi \circ \iota \right\| <\delta.
\]
 From this we can conclude (see \cite [Lemma 2.3]{EH}) that
\[
\left\| \iota-\psi \circ \iota \right\|_{cb} < \delta n.
\]
We next perturb $\psi $ in order to satisfy
\begin{equation}
  \psi |_{F} = id_{F}.   \label {E5}
\end{equation}
It follows from \cite [Lemma 5.2]{ER4}  that there is a projection
$P$ of $\tilde F$ onto $F$ with
$1\leq \left\| P\right\| _{cb}\le n^{2}$.
Then
\[
\tilde \psi =(\iota-\psi) \circ P+\psi: \tilde F \rightarrow \A
\]
is a completely bounded mapping satisfying (\ref{E5}).
Given  any $v\in V$,  $P\circ \al(v)$ is an element in
$F = \al(V) \cap \A$.
Then there exists $v_0\in V$ such that $P\circ \al (v) = \al (v_0)$.
Since $\psi$ satisfies (\ref {E4}), we have
\[
\tilde \be \circ  \psi \circ P \circ \al (v) =
\tilde \be \circ  \psi \circ \al (v_0)   = \tilde \be \circ \al (v_0)
=\tilde \be \circ  P \circ \al (v) ,
\]
and thus
\[
\tilde \be \circ \tilde \psi  \circ \al =
\tilde \be \circ (\iota -\psi)\circ P\circ \al
+\tilde \be\circ\psi\circ \al
= \tilde \be \circ \psi \circ \al = \varphi.
\]
This shows that
$\tilde \psi$ also satisfies (\ref {E4}).
Finally, $\tilde \psi$  satisfies
\[
\| \tilde \psi \| _{cb}
\le \left\| \iota -\psi \circ \iota\right\|_{cb}
\left\| P\right\| _{cb}+(1+\delta )
\le \delta n^{3}+(1+\delta)<1+\e .
\]
If we let  $\tilde \al = \tilde \psi \circ \al $ and let $\tilde \be$
also denote
its restriction to $\A$, we obtain the
completely bounded factorization
\[
\begin{array}{ccccc}
&  & {\A} &  &  \\
& {\scriptstyle \tilde \al }\longnearrow  &  & \longsearrow
{\scriptstyle
\tilde \be } &  \\
{V} \hspace{-10 pt} &  & \stackrel{\varphi}{\hspace{-10pt}{\hbox to 30pt
{\rightarrowfill}}\hspace{-10pt} } &  &\hspace{-10 pt} W
\end{array}
\]
such that
\[
\|\tilde \al\|_{cb} \le \|\tilde \psi \|_{cb}\|\al\|_{cb} <
(1+\e)\|\al\|_{cb}
\]
and
\[
\|\tilde \be\|_{cb}< (1+\e)\|\be\|_{cb}.
\]
We have $\al(V) \cap \A \subseteq \tilde \al (V)$
since for any $x \in F =  \al(V)\cap {\A}$,
there exists $v \in V$ such that
\[
x = \al (v) = \tilde \psi \circ \al(v) = \tilde \al (v) \in \tilde
\al(V),
\]
where we used the fact that $\tilde \psi |_{F} = id_{F}$.
\endproof

\begin{theorem}
\label{PG.3}
Let $\A$ be a locally reflexive $C^*$-algebra and let $E$ be a
finite-dimensional
subspace  of  $\A$.
If $\tilde F$ is a finite-dimensional subspace of $\A^{**}$ such that
$E\subseteq \tilde F$ and
\[
d_{cb}(\tilde F, \B) < \lambda
\]
for some finite-dimensional $C^*$-algebra $\B$, then there exists a
finite-dimensional subspace  $ F$ of $\A$ such that
$E \subseteq F$ and
\[
d_{cb}(F, \B) < \lambda.
\]
\end{theorem}
\proof
Let $\B$ be a finite-dimensional $C^*$-algebra and let
\[
\al : \B \to  \tilde F \subseteq \A^{**}
\]
be a  linear isomorphism from $\B$ onto $\tilde F$ such that
$\|\al\|_{cb} \|\al^{-1}\|_{cb} < \la$.
Since $\B$ is an injective operator space, the inverse mapping
$\al^{-1} : \tilde F \to \B$ has a completely bounded
extension $\be : \A^{**} \to \B$ with $\|\be\|_{cb} =
\|\al^{-1}\|_{cb}$.
Then we obtain a completely bounded factorization
\[
\begin{array}{ccccc}
&  & \A^{**} &  &  \\
& {\scriptstyle \al}\longnearrow  &  & \longsearrow {\scriptstyle
\be} &  \\
\B\hspace{-10 pt} &  & \stackrel{id_{\B}}{\hspace{-10pt}{\hbox to 30pt
{\rightarrowfill}}\hspace{-10pt} } &  &\hspace{-10 pt} \B
\end{array}
\]
with $\|\al\|_{cb} \|\be\|_{cb} < \lambda$.
It follows from Lemma \ref {PG.2} that
we may find a completely bounded factorization
\[
\begin{array}{ccccc}
&  & \A &  &  \\
& {\scriptstyle \tilde \al}\longnearrow  &  & \longsearrow {\scriptstyle
\tilde \be} &  \\
\B\hspace{-10 pt} &  & \stackrel{id_{\B}}{\hspace{-10pt}{\hbox to 30pt
{\rightarrowfill}}\hspace{-10pt} } &  &\hspace{-10 pt} \B
\end{array}
\]
such that  $\|\tilde \al\|_{cb} \|\tilde \be\|_{cb} < \lambda$
and
\[
E \subseteq \tilde F \cap \A = \al(\B) \cap \A \subseteq \tilde \al(\B).
\]
Then $ F = \tilde \al(\B)$ is a
finite-dimensional subspace of $\A$ containing $E$ and
$\tilde \al  : \B \to F$ is a linear isomorphism with
$\tilde \al^{-1} =   \tilde\be | _F$.
 From this we conclude that
\[
d_{cb}(F, \B)  \le \|\tilde \al \|_{cb} \|\tilde \al^{-1} \|_{cb} < \lambda.
\]
\endproof

The following result is already known to Blackadar and Kirchberg \cite {BK2}
(see our discussion in $\S 2$).
For the convenience of the readers, we include the following simple proof.
Here, we do not need to assume the separability for $\A$.

\begin{theorem}
\label{PG.0}
If $\A$ is a subhomogeneous $C^*$-algebra, then
$\A$ is an $\OL_{\infty, 1^+}$ space.
\end{theorem}
\proof
If $\A$ is a subhomogeneous $C^*$-algebra, then
all irreducible  representations of $\A$  are finite-dimensional
with ${\rm dim} \le n$ for some  positive integer $n\in \N$, and thus
its second dual  $\A^{**}$ has the form
\[
\A^{**} \cong \oplus _{k=1}^n L_\infty(X_k, {\mathcal M}_k, \mu_k)
\check \otimes M_{k},
\]
where  we simply assume that $L_\infty(X_k, {\mathcal M}_k, \mu_k) = \{0\}$
if $\A$ does not have any irreducible representation of dimension $k$.
Since unital commutative $C^*$-algebras are rigid $\OL_{\infty, 1^+}$
spaces (see
discussion in $\S 1$),  it is easy to see that
$\A^{**} \cong \oplus _{k=1}^n L_\infty(X_k, {\mathcal M}_k, \mu_k)
\check \otimes M_{k}$ is a rigid ${\OL}_{\infty, 1^+}$ space (or a
strong NF algebra).
The $C^*$-algebra  $\A$ is nuclear and thus is locally reflexive.
Given any finite-dimensional subspace $E\subseteq \A$,
which can be  regarded as a finite-dimensional subspace of $\A^{**}$, and
any $\e> 0$, there exists a finite-dimensional subspace $\tilde F$
in $\A^{**}$ such that $E \subseteq \tilde F$ and
\[
d_{cb}(\tilde F, \B) < 1+\e
\]
for some finite-dimensional $C^*$-algebra $\B$.
It follows from Theorem \ref {PG.3} that we may find a finite-dimensional
subspace $F$ in $\A$ such that $E \subseteq F$ and
\[
d_{cb}(F, \B) < 1+\e.
\]
This shows that $\A$ is an $\OL_{\infty, 1^+}$ space.
\endproof

\section{Non-subhomogeneous nuclear $C^*$-algebras}

\begin{lemma}
\label {PNG.1}
If $\A$ is a non-subhomogeneous $C^*$-algebra, then there is a
completely isometric  and completely order preserving injection
\[
\theta : {\B}(\ell_2) \to \A^{**},
\]
which is a weak$^*$ homeomorphism from ${\B}(\ell_2)$ onto
$\theta({\B}(\ell_2))$.
\end{lemma}
\label{PNG.2}
\proof
If $\A$  is non-subhomogeneous, then either $\A$ has an
infinite-dimensional irreducible representation, or
all irreducible representations of $\A$ are finite-dimensional,
but the dimensions are not uniformly bounded.

If $\A$ has an infinite-dimensional irreducible representation
$\pi : \A \to {\B}(H)$ with ${\rm dim} H = \infty$,
then $\pi$ induces a  unique  normal (i.e. weak$^*$-continuous)
representation
$\tilde \pi: \A^{**}\to \B(H)$ from $\A^{**}$ onto $\B(H) = \pi(\A)''$,
and there is a  central projection $p\in \A^{**}$
such that ${\rm ker}\tilde \pi = (1-p)\A^{**}$,
where ${\rm ker} \tilde \pi$ is the   kernel of $\tilde \pi$.
This induces a normal $*$-isomorphism
\[
\B(H) \cong p \A^{**}
\]
(see Takesaki's book \cite {Takesaki}).
Since ${\rm dim} H = \infty$ (which could be uncountable), we may
identify $\ell_2$ with a subspace of $H$ and identify $\B(\ell_2)$ with
a von Neumann  subalgebra of $\B(H)$.
In this case, we obtain an injective  normal  $*$-homomorphism
\[
\theta:  \B(\ell_2) \to \A^{**}.
\]
It is obvious that $\theta$ is a completely isometric and completely order
preserving weak$^*$ homeomorphism from $\B(\ell_2)$ onto $\theta(\B(\ell_2))$.

Now if we assume that all irreducible representations
of $\A$  are finite-dimensional, but not uniformly bounded,
then there exists a strictly
increasing sequence of positive integers $n(k) \in \Bbb N$
such that for each $k$, $\A$ has an irreducible representation
\[
\pi^k : \A \to M_{n(k)}.
\]
By the same reason as above, for each $k \in \N$ there exists a
central projection $p^k$ in $\A^{**}$
such that
\[
M_{n(k)} \cong p^k \A^{**}.
\]
Since $\pi^k$ are distinct irreducible representations of $\A$ and
there is no  non-trivial central projections in $M_{n(k)}$,
$\{p^k\}$ must be all distinct and  mutually orthogonal.
Then $p = \sum_{k=1}^\infty p^k$ is a central projection in $\A^{**}$ and
$\prod_{k=1}^\infty  M_{n(k)}$ can be identified with the von Neumann
subalgebra
$p\A^{**}$ in $\A^{**}$.
For each $k \in \Bbb N$, we let $P^k$ denote the natural truncation
mapping
\[
P^k: x \in {\B}(\ell_2) \to x^{n(k)}\in M_{n(k)}.
\]
It is known (see \cite {ER1}) that  the canonical mapping
\[
\theta : x \in {\B}(\ell_2) \to \theta(x) = (P^k(x)) \in
\prod_{k=1}^\infty  M_{n(k)}
\cong p\A^{**}
\]
is a completely positive and completely isometric injection.
The mapping $\theta$ is  a complete order isomorphism
from $\B(\ell_2)$ onto $\theta(\B(\ell_2))$ since for every
$x \in \B(\ell_2)$, $x\ge 0$ if and only if
$P^k(x) \ge 0$ for all $k \in \N$.
It is weak$^*$ continuous since each truncation mapping
$P^k$ is weak$^*$-continuous.
It is also easy to see that $\theta(\B(\ell_2))$ is weak$^*$ closed in
$\prod_{k=1}^\infty  M_{n(k)} \cong p\A^{**}$.
Then $\theta$ is a weak$^*$ homeomorphism from
$\B(\ell_2)$ onto $\theta(\B(\ell_2))$.
\endproof

In the following,  we let $\tau (x) = \int_0^1 x(t) dt$ denote the normal
(tracial) state on $L_\infty [0, 1]$
and let  $\langle x ~ | ~ y \rangle = \tau (y^* x) $ stand for the
inner product on $L_2 [0, 1]$.
We let  $r_k(t)={\rm sgn} \sin(2^k\pi t)$ denote the Rademacher functions
on $[0, 1]$.
Then $\{r_k\}$ is a sequence of  self-adjoint unitary
elements in $L_\infty [0, 1]$ since each  $r_k$ is a  real
valued function on $[0, 1]$
with $r_k^2 = 1$.
Moreover, $\{r_k\}$ forms an orthonormal set in $L_2 [0, 1]$
since $\langle r_k ~ | ~ r_{k'}\rangle = \delta_{k, k'}$.

\begin{lemma}
\label{PNG.2}
%\label{rad}
For any $x\in L_\infty [0,1]$, we have
\[
\lim_{k\to \infty}  \tau(r_k x)= 0.
\]
\end{lemma}
\proof
Since $\{r_k\}$ is an orthonormal set in $L_2[0, 1]$, we have from
the Bessel inequality
that for any $x \in L_\infty [0, 1]$
\[
\sum_{k =1}^\infty |\tau(r_k x)|^2 =
\sum_{k=1}^\infty |\langle x ~ | ~ r_k \rangle |^2 \le \|x\|^2_2 .
\]
This implies $\lim_{k \to \infty}  \tau (r_k x)= 0$.
\endproof

\begin{lemma}
\label {PNG.3}
%\label{2}
Let $F \subseteq L_\infty[0,1]\ten M_n$
and $S \subseteq (L_\infty[0,1]\ten M_n)_*$ be finite-dimensional
operator spaces.  Then for  every $\e >0$,
there exists a complete contraction
$u: L_\infty[0,1]\ten M_n\to M_n$ and  a complete isometry
$v: M_n\to L_\infty[0,1]\ten M_n$ such that $u\circ v=id_{M_n}$  and
\begin{enumerate}
  %\item[(i)] $\|u\|_{cb}\le 1$ and  $\|v\|_{cb}\le 1$,
  \item[(i)]  $\|u(x)\|_{M_n} <  \e  \| x\|$,
  \item[(ii)] $| \langle s,~   v(y)\rangle| < \e \|s\| \| y\|$
   \end{enumerate}
for all $x\in F$, $y\in M_n$ and $s\in S$.
\end{lemma}
\proof  Using the Rademacher functions,
we may define a sequence of complete contractions
$u_k :  L_\infty[0,1]\ten M_n \to M_n$ by letting
\[
u_k(x) = (\tau \otimes id)((r_k\otimes 1)x)
\]
for all $x \in L_\infty[0,1]\ten M_n$.
Let $\{e_{ij}\}$ denote the matrix unit of $M_n$.
For every  $x \in  L_\infty[0,1]\ten M_n$,  we can
write
\[
x = \sum_{i,j}^n x_{ij} \otimes e_{ij}
\]
with $x_{ij} \in L_\infty([0,1])$.
It follows from Lemma \ref {PNG.2} that
\begin{equation}
\label {F.point-norm}
u_k(x) = \sum_{i,j}^n \tau(r_k x_{ij})  e_{ij} \to 0
\end{equation}
in $M_n$.
Since $F$ is a finite-dimensional subspace  of $L_\infty[0,1]\ten M_n$,
its closed unit ball $D_F$  is norm compact and we can conclude
from (\ref {F.point-norm}) that  $u_k \to 0$ uniformly on $D_F$.
Then for every $\varepsilon > 0$,
there exists $k_1 \in \N$ such that for all  $k \ge k_1$,
$u_k$ satisfy the condition (i), i.e.
\[
\|u_k(x)\|_{M_n} \le \e \|x\|
\]
for all $x \in F$.

We can consider another sequence of mappings
$v_k: M_n \to  L_\infty[0,1]\ten M_n$ given by
\[
v_k(y) = r_k \otimes y
\]
for all $y \in M_n$.
Since $r_k^2 = 1$,
it is clear that $v_k$ are completely isometric injections,
and  $u_k \circ v_k = id_{M_n}$ since
\[
u_k\circ v_k(y) = (\tau\otimes id) (r_k^2 \otimes  y) = \tau(1)y  = y
\]
for all $y \in M_n$.

To see (ii),  let us first recall from \cite {ER1} that there is a
complete isometry
\[
(L_\infty[0,1]\ten M_n)_* \cong L_1[0, 1]\hat \otimes T_n,
\]
where $L_1[0, 1]\hat \otimes T_n$ is the operator projective tensor
product of
$L_1[0, 1]$ and $T_n $ introduced in \cite {BP} and \cite {ER2}.
Given any $s \in S\subseteq (L_\infty[0,1]\ten M_n)_* $, we may write
\[
s =  \sum_{i,j}^n s_{ij} \otimes   \tilde e_{ij},
\]
where $\{\tilde e_{ij}\}$ is the canonical dual basis of $\{e_{ij}\}$ in
$T_n$
and $s_{ij}$ are integrable functions in $L_1[0, 1]$ with
$\|s_{ij}\|_{L_1[0, 1]} \le \|s\|$.
Since $L_\infty [0, 1]$ is norm dense in $L_1[0, 1]$,
there exist $h_{ij}\in L_\infty [0,1]$ such that
\[  \|s_{ij}- h_{ij}\|_{L_1[0, 1]} <  {\frac {\e}{4n^2}} \|s\|.
\]
It follows that
\begin{equation}
\label {FNG.1}
|\tau(r_k s_{ij} - r_k h_{ij})| \le \|r_k\|_{L_{\infty}[0,1]}
\|s_{ij} - h_{ij}\|_{L_1[0,1]} < {\frac {\e}{4n^2}} \|s\|.
\end{equation}
Since  $\lim_{k\to \infty} |\tau (r_k h_{ij})| = 0$, we may
choose $k_2 \ge k_1$ such that
\begin{equation}
\label{FNG.2}
|\tau(r_{k_2} h_{ij})|< {\frac {\e}{4n^2}}\|s\| .
\end{equation}
Then given any  $y = \sum_{i, j}^n y_{ij} \otimes e_{ij} \in {M_n}$, we
deduce  from (\ref {FNG.1}) and (\ref {FNG.2}) that
\begin{eqnarray*}
|\langle s, ~ v_{k_2}(y)\rangle|  &=& |\sum_{i,j}^n \tau (r_{k_2} s_{ij})y_{ij}|
\\
&\le&   \sum_{i,j}^n |y_{ij}|\, |\tau (r_{k_2} s_{ij} - r_{k_2} h_{ij})|
   + \sum_{i,j}^n |y_{ij}|\,  |\tau (r_{k_2} h_{ij})| \\
&<  & {\frac \e 2} \|s\| \, \|y\|.
\end{eqnarray*}
Since the closed unit sphere $D_S$  of $S$  is norm compact,
we can conclude that
\[| \langle s,~ v_{k_2}(y)\rangle| < \e \|s\| \| y\|
\]
for all $y\in M_n$ and $s\in S$.
Then the mappings $u = u_{k_2}$ and $v = v_{k_2}$ satisfy the
conditions (i) and (ii).
\endproof

As we discussed in $\S 1$, a $C^*$-algebra  $\A$ is  nuclear
if and only if there exist approximate diagrams of complete contractions
in (\ref {F1.nuc}) which approximately commute in
the point-norm topology.
This is equivalent to saying that for any finite-dimensional
subspace $E$  of $\A$ and
$\e > 0$, there exists a matrix space $M_n$ and
a commutative diagram of completely bounded mappings
\begin{equation}
\label {FNG.3}
\begin{array}{ccccc}
&  & M _{n} &  &  \\
& {\al }\longnearrow  &  & \longsearrow {\scriptstyle
\be} &  \\
E\hspace{-10 pt} &  & \stackrel{\iota}{\hspace{-10pt}{\hbox to 30pt
{\rightarrowfill}}\hspace{-10pt} } &  &\hspace{-10 pt} \A
\end{array}
\end{equation}
such that $\|\al\|_{cb} \le 1$ and $\|\be\|_{cb} < 1+\e$.

The following Lemma is due to Oikhberg together with the
observation that the modification  of the conclusion
is true for operator spaces
with completely bounded approximation property, but
fails for general operator spaces.
For completeness, we include the proof.

\begin{lemma}
\label{Oik} Let $\A$ be a nuclear $C^*$-algebra, then for every
finite-dimensional subspace $F\subseteq \A$ and $\e>0$, there exists a
finite-codimensional subspace $W$ of $V$ such that the quotient
mapping $q: \A \to \A/W$ induces a completely contractive linear
isomorphism
\[
q|_F : F \to (F+W)/W
\]
with $\| q|_F^{-1} \|_{cb} < (1+\e)$.
\end{lemma}
\proof
Since $\A$  is nuclear,
it follows from (\ref {FNG.3}) that  for every finite-dimensional
subspace
$F \subseteq \A$, there is a finite-rank mapping
\[
T : \A\to \A
\]
such that $\|T\|_{cb} < (1+\e)$ and $T(x)=x$ for all $x\in F$.
Then  $W={\rm ker}\,T$  is a finite-codimensional subspace of $\A$,
and $T$ determines a canonical mapping $\hat T: \A/W\to \A$
given by
\[
\hat T(x+W) = T(x).
\]
Since $q$ is a complete quotient mapping from $\A$ onto $\A/W$, we have
\[
\|\hat T\|_{cb} = \|T\|_{cb} < (1+\e),
\]
and
\[
\hat{T}\circ q(x) = \hat T(x+W)= T(x) = x
\]
for all $x \in F$.
This shows that $\hat T$ restricted to $q (F) = (F+W)/W$ is a left
inverse of
$q|_F$.
Then $q |_F : F \to (F+W)/W$ is a completely contractive linear
isomorphism with
\[
\| q  |_F^{-1}\|_{cb} \le \|\hat T \|_{cb} < (1+\e).
\]
\endproof

\begin{theorem}
\label{PNG.4}
Let $\A$ be a non-subhomogeneous nuclear $C^*$-algebra.
Then for any finite-dimensional subspace
$E\subseteq A$ and  $0<\e<\frac 12 $,
there exists a subspace $\tilde F \subseteq \A^{**}$
containing $E$ and a linear isomorphism
\[
\p: M_n \to \tilde F
\]
such that
\[
\|\p\|_{cb} < 3 +2\e ~ \mbox{and } ~ \|\p^{-1}\|_{cb} < 2+ 12\e.
\]
\end{theorem}
\proof
Let $E$ be a  finite-dimensional subspace of $\A$.
Since $\A$ is nuclear,  it follows from (\ref{FNG.3}) that
there exists a matrix space $M_n$ and
a commutative diagram of completely bounded mappings
\[
\begin{array}{ccccc}
&  & M _{n} &  &  \\
& {\al }\longnearrow  &  & \longsearrow {\scriptstyle
\be} &  \\
E\hspace{-10 pt} &  & \stackrel{\iota}{\hspace{-10pt}{\hbox to 30pt
{\rightarrowfill}}\hspace{-10pt} } &  &\hspace{-10 pt} \A
\end{array}
\]
such that $\|\al\|_{cb} \le 1$ and $\|\be\|_{cb} < 1+\e$.
For the  finite-dimensional
subspace  $F=\be(M_n)$ of $\A$, there exists a finite-codimensional
subspace $W$ of $\A$ such that
the quotient mapping  $q: \A \to A/W$ induces a
completely contractive linear isomorphism
\[
q|_F : F \to (F+W)/W
\]
with  $\|q|_F^{-1}\|_{cb} < 1+\e$.
In this case,
\[
S = W^\perp \cong (A/W)^*
\]
is a finite-dimensional subspace of $\A^*$.

Since we may identify $ L_\infty [0, 1] \ten M_n$ with a von Neumann
subalgebra of $\B(\ell_2)$, the canonical mapping
$\theta: \B(\ell_2)\to \A^{**}$ in
Lemma \ref {PNG.1} induces a  weak$^*$ continuous
completely isometric injection from $L_\infty [0, 1] \ten
M_n$ into $\A^{**}$, which is still denoted by $\theta$.
Since $ L_\infty [0, 1] \ten M_n$ is
an injective operator space, there is a complete contraction
\[
P:\A^{**}\to L_\infty [0, 1] \ten M_n
\]
such that $P\circ \theta = id_{L_\infty [0, 1] \ten M_n}$.
Applying  Lemma \ref {PNG.3} to the finite-dimensional spaces
$P(F)$ in $L_\infty [0, 1] \ten M_n$ and $\theta_*(S)$ in
$(L_\infty [0, 1] \ten M_n)_*$, and $\e'=\frac{\e}{n^2}$,
we may find a complete contraction
$u: L_\infty [0, 1] \ten M_n \to M_n$ and a complete isometry
$v: M_n\to L_\infty [0, 1] \ten M_n$
such that  $u \circ v = id_{M_n}$ and
\begin{enumerate}
  \item[(i)]  $\|u\circ P(x)\| <  \e'  \| P(x)\| \le \e' \|x\|$,
  \item[(ii)] $| \langle  s,~   \theta\circ v(y)\rangle| <
  |\langle  \theta_*(s),~   v(y)\rangle| < \e' \|\theta_*(s)\| \|y\| \le
\e' \|s\| \|y\|$
  \end{enumerate}
for all $ x\in F$, $y\in M_n$ and $s\in S$.

The mapping  $T = u \circ  P : \A^{**} \to M_n$  is a complete
contraction such that
$\|T|_F\| < \e'$ by (i).
Since ${\rm dim} F \le n^2$, it follows from [Lemma 2.3]\cite {EH}  that
\[
\|T|_F\|_{cb} \le n^2 \|T |_F\| < \e,
\]
and thus
\begin{equation}
\label {EstT1}
\|T_m(x)\| \le \|T|_F\|_{cb} \|x\| < \e \|x\|
\end{equation}
for all  $x \in M_m(F)$ and $m \in \N$.
If we let  $G = \theta\circ v (M_n)$,
then $G$ is a finite-dimensional subspace of $\A^{**}$,
and
\[
\theta\circ v : M_n \to G
\]
is a completely isometric isomorphism.
Given any $\tilde y \in M_m(G)$, there exists a unique $y \in M_m(M_n)$
such that
$\tilde y = \theta_m \circ v_m(y)$, and
\begin{equation}
\label{EstT2}
\|T_m(\tilde y)\| = \|u_m \circ P_m\circ \theta_m \circ v_m(y) \|
= \|u_m \circ v_m(y) \| = \|y\|  =\|\tilde y\|.
\end{equation}

On the other hand, we let
\[
R= q^{**} : \A^{**} \to (\A/W)^{**} \cong \A/W
\]
denote the second adjoint of  $q$.
Since $q^* : (\A/W)^* \to W^\perp \subseteq \A^*$ is a complete
isometry from  $(\A/W)^*$ onto $W^\perp$, we have from (ii) that  for
any $\tilde y = \theta\circ v(y) \in G$ with $y \in M_n$,
\begin{eqnarray*}
\|R(\tilde y)\| &=& \sup \{|\langle R(\tilde y), ~ s \rangle|:
~ s \in (\A/W)^*, \|s\| \le 1\}\\
&=& \sup \{|\langle \tilde y, ~ q^*(s) \rangle|:
~ s \in (\A/W)^*, \|s\| \le 1\}\\
&=& \{|\langle \theta \circ v (y), ~ \tilde s) \rangle|:
~ \tilde s = q^*(s) \in S = W^\perp, \|\tilde s\| \le 1\}\\
%&=& \{|\langle  v (y), ~ \theta_*(\tilde s) \rangle|:
%~ s \in S = W^\perp, \|s\| \le 1\} \\
&<&  \e' \|y\| = \e' \|\tilde y\|.
\end{eqnarray*}
This shows that
\[
\|R|_G\|_{cb} \le n^2 \|R|_G\| < \e
\]
and thus
\begin{equation}
\label {EstR1}
\|R_m(\tilde y)\| \le \|R|_G\|_{cb}\, \|\tilde y\| < \e \|\tilde y\|
\end{equation}
for all $\tilde y  \in M_m(G)$.
Given any $x \in M_m(F)$, we also have
\begin{equation}
\label {EstR2}
\|x\| = \|(q|_F^{-1})_m \circ (q|_F)_m(x)\| < (1+\e) \|(q|_F)_m(x)\|
= (1+\e)  \|R_m(x)\|.
\end{equation}

Now we are ready to define the mapping
\[
\p = \be + \theta\circ v\circ (id_{M_n} - \al \circ \be): M_n \to
\A^{**}.
\]
It is clear that $\p$ is completely bounded with
\begin{eqnarray*}
\|\p\|_{cb} &\le&   \|\be\|_{cb} +
\|\theta\circ v\circ (id_{M_n} - \al \circ \be)\|_{cb} \\
&=&  \|\be\|_{cb} +
\|id_{M_n} - \al \circ \be\|_{cb}  \\
&< & (1+\e) + 1 + (1+\e) = 3 + 2\e.
\end{eqnarray*}
We let $\tilde F = \p(M_n)$.
If $x \in E$, then $\al(x) \in M_n$ and
\[
\p(\al(x)) = \be(\al(x)) + \theta\circ v\circ (id_{M_n} - \al \circ \be)
(\al(x))
= x.
\]
This shows that $E = \p(\al(E)) \subseteq \tilde F$.

The mapping  $\p$ is a linear isomorphism from $M_n$ onto $\tilde F$.
To see this, we  first claim  that $F \cap G = \{0\}$.
Given any  $x \in  F\cap G$,
we have from (\ref{EstR1}) and (\ref{EstR2}) that
\[
\|x\| < (1+\e) \| R(x)\|
  <  \e(1+\e) \| x \|  < 2\e \| x \|.
\]
Since $\e<\frac 12$, we must have $\|x\|=0$ and thus $F\bigcap G =
\{0\}$.
If $y \in M_n$ such that $\p(y) = 0$, then
\[
\theta\circ v\circ (id_{M_n} - \al \circ \be) (y) = -\be(y) \in F \cap G
= \{0\}.
\]
This implies that $y = \al \circ \be (y) = 0$.
Therefore,  $\p$ is an injection and thus a linear isomorphism from
$M_n$ onto $\tilde F$.

Finally, we want  to show that $\|\p^{-1}\|_{cb} < 2 + 12\e$.
Given any  $x \in M_m(F)$ and $y \in M_m(G)$,
we have from (\ref {EstR1}) and (\ref {EstR2}) that
\begin{eqnarray*}
  \frac{1}{1+\e}   \| x\| &< &  \| R_m(x)\|
   \le  \| R_m(x+y)\|+ \| R_m(y)\| \\
    &< & \| x+ y \| + \e \| y\| .
\end{eqnarray*}
Similarly, we have
\begin{eqnarray*}
   \|  y\| &= &  \| T_m(y)\|
   \le  \| T_m(x+y)\|+ \| T_m(x)\| \\
    &<& \| x+y\| + \e \| x\| .
\end{eqnarray*}
 From this we can conclude that
\[
\| x\| < (1+\e) \| x+y\|  + \e (1+\e)\|x+y\| + \e^2(1+\e)\|x\|,
\]
and thus
\[
\| x\| < \frac{(1+\e)^2 }{1-\e^2(1+\e)} \| x+y\| .
\]
Similarly, we have
\[
  \| y\| < \| x+y\| + \e(1+\e) \| x+y \| +\e^2(1+\e) \| y\|,
\]
and thus,
\[
\|y\| < \frac{(1+\e)^2  }{1-\e^2(1+\e) } \| x+y \| .
\]
Since $\e < \frac 12$, we have
\[
{\frac{(1+\e)^2  }{1-\e^2(1+\e) }} = 1 + {\frac{(1+\e)^2  - 1 +
\e^2(1+\e)}{1-\e^2(1+\e) }}
= 1 + \frac{2\e + 2\e^2 + \e^3}{1-\e^2(1+\e) } < 1 +  6\e,
\]
and thus
\begin{equation}
\label{FNG.max}
\mbox{max}\{\|x\|, \|y\|\} < (1+6\e) \|x+y\|.
\end{equation}
Since $\p(M_n)\subset F+G$,  for any  $y\in M_m(M_n)$ we have
\begin{eqnarray*}
  \| y\| &=& \| \al_m \circ \be_m (y) + y-\al_m \circ \be_m  (y)\| \\
&\le& \|\al_m\circ\be_m (y) \|+\|u_m \circ v_m ( y-\al_m \circ\be_m
(y))\| \\
   &\le& \| \be_m (y)\| + \| v_m ( y-\al_m  \circ \be_m  (y))\|  \\
&\le& \| \be_m (y)\|+\|\theta_m \circ v_m ( y-\al_m\circ \be_m (y))\|
\\
   &\le& 2 (1+ 6\e)
\| \be_m (y)+ \theta_m \circ v_m ( y-\al_m  \circ \be_m  (y))\|  \\
&=&  2  (1+6\e)  \| \p_m(y)\|.
\end{eqnarray*}
This shows that
\[
\|\p^{-1}\|_{cb} \le 2 + 12\e.
\]
\endproof

\begin{theorem}
\label{PNG.6}
If $\A$ is a non-subhomogeneous nuclear $C^*$-algebra, then
$\A$ is an $\OL_{\infty, \la}$ space for every $\la > 6$.
\end{theorem}
\proof
Given any  $\la > 6$, we may find a sufficiently small
$\e$ such that $\e < {\frac 12}$
and $(3+2\e) (2+ 12\e) < \la$.
For  any finite-dimensional subspace $E \subseteq
\A$, it follows from Theorem
\ref {PNG.4} that  there exists a
finite-dimensional subspace $\tilde F \subseteq \A^{**}$ containing
$E$ and a linear isomorphism
$\p: M_n\to \tilde F$ with
$\|\p\|_{cb} < 3 + 2\e$ and $\|\p^{-1}\|_{cb} < 2+12\e$.
Then we have
\[
d_{cb}(\tilde F, M_n) \le \|\p\|_{cb} \|\p^{-1}\|_{cb} < (3 +
2\e)(2+12\e) < \la.
\]
It follows from
Theorem \ref {PG.3} that there exists a finite-dimensional
subspace $F$ of $\A$ such that $E \subseteq F$ and
\[
d_{cb}(F, M_n) < \la.
\]
This shows that  $\A$ is an ${\OL}_{\infty, \la}$ space.
\endproof

Combining Theorem \ref {PG.0} and Theorem \ref {PNG.6}, we can
conclude that every $C^*$-algebra $\A$ (either subhomogeneous or
non-subhomogeneous) is an $\OL_{\infty, \la}$ space for every $\la > 6$.
This completes the proof for Theorem \ref {P1.2}.

\section{remarks and questions}

As we discussed in $\S 1$,  we can define an
invariant
\[
{\OL}_\infty(\A) = \mbox{inf}\{\la > 1: ~ \A ~ \mbox{is ~ an~ }
{\OL}_{\infty, \la} ~ \mbox{space}\}
\]
for every nuclear $C^*$-algebra $\A$.
In general, we have $1 \le {\OL}_\infty(\A) \le 6$, and
${\OL}_\infty(\A) = 1$ if and only if $\A$ is an $\OL_{\infty, 1^+}$
space.
This includes the case when $\A$ is a rigid $\OL_{\infty, 1^+}$ space
or a strong NF algebra.

\begin{question}
\label {Q6.1}
It would be interesting to know whether $\OL_{\infty, 1^+}$  implies rigid
$\OL_{\infty, 1^+}$ on unital nuclear $C^*$-algebras.
\end{question}

It is known from Theorem \ref {P2.quasi} that if $\A$ is a
unital $C^*$-algebra with $\OL_{\infty}(\A) =1$,
then $\A$ is nuclear and quasi-diagonal.

\begin{question}
\label {Q6.2}
Does a nuclear and quasi-diagonal unital $C^*$-algebra must be
an $\OL_{\infty, 1^+}$ space ?
\end{question}

To investigate Questions \ref {Q6.1} and \ref {Q6.2},
one could look at the $C^*$-algebra $\B n$  introduced by L. Brown
\cite {Br}, which is an essential extension
\[
0\to K_\infty \to \B n \to C({\Bbb R}P^2) \to 0
\]
of the continuous functions on the real projective plane ${\Bbb R}P^2$
by $K_\infty = K(\ell_2)$.
One could also look at the $C^*$-algebra $C^*(s \oplus s^*)$ generated by the
direct sum of the unilateral shift $s$ (on $\ell_2$) and its adjoint $s^*$,
for which we have the extension
\[
0\to K_\infty \oplus K_\infty \to C^*(s\oplus s^*) \to C(S^1) \to 0.
\]
It is known (see Blackadar and Kirchberg \cite [$\S 2$]{BK2})  that
$\B n$ and $C^*(s \oplus s^*)$ are nuclear quasi-diagonal unital
$C^*$-algebras, but they are not inner quasi-diagonal and thus are not
rigid $\OL_{\infty, 1^+}$ spaces.
It would be interesting to calculate $\OL_\infty (\B n)$ and
$\OL_\infty (C^*(s \oplus s^*))$.
If $\OL_\infty(\B n) = 1$ (or $\OL_\infty(C^*(s \oplus s^*)) = 1$),
we would get an example of a unital $C^*$-algebra, which is an
$\OL_{\infty, 1^+}$ space,  but is not a rigid  $\OL_{\infty, 1^+}$ space.
On the other hand, if $\OL_\infty(\B n) > 1$ (or
$\OL_\infty(C^*(s \oplus s^*)) > 1$),  then we would obtain an example of a
unital $C^*$-algebra, which is  nuclear and quasi-diagonal, but is not
$\OL_{\infty, 1^+}$.
This investigation would either give us  a negative answer to Question \ref
{Q6.1} or give us a negative answer to  Question \ref {Q6.2}.

On the other hand, it is known from Theorem \ref {P.stablef} that
if $\A$ is a  non-stably finite nuclear unital $C^*$-algebra,
then we have
\[
(\frac {1+ \sqrt 5} 2)^{\frac 1 2} < {\OL}_\infty(\A) \le 6.
\]
For example, we may  consider the Cuntz algebra ${\mathcal O}_n$
(with $2 \le n < \infty$), or the  Toeplitz algebra ${\mathcal T}(S^1)$
on the unit ball of  ${\Bbb C}$.
Since ${\mathcal T}(S^1)$ is the $C^*$-algebra generated by the
unilateral shift $s$ (on $\ell_2$) and has an extension
\[
0 \to K_\infty \to {\mathcal T}(S^1) \to C({S^1}) \to 0,
\]
it is an infinite nuclear unital $C^*$-algebra, and thus  we have
\[
(\frac {1+ \sqrt 5} 2)^{\frac 1 2} < {\OL}_\infty({\mathcal T}(S^1)) \le 6.
\]
This example shows that the (rigid) ${\OL}_{\infty, 1^+}$ structure is not
preserved by $C^*$-algebra extensions.
We may also consider the  Toeplitz algebra ${\mathcal T}(S^3)$ on the
unit ball of ${\Bbb C}^2$, which is a finite, but not 2-finite
nuclear unital $C^*$-algebra (see \cite [$\S 6.10$]{Bla}).

\begin{question}
It would be interesting to know if there is any stably finite unital
$C^*$-algebra for which $\OL_\infty (\A) >1$, or
$\OL_\infty (\A)> (\frac {1+ \sqrt 5} 2)^{\frac 1 2}$.
\end{question}

Finally we remark that in a recent work of Junge, Nielsen, Ruan and Xu
\cite {four},
we  obtained the constant $\OL_{\infty}(\A) \le 3$ for every nuclear
$C^*$-algebra  with completely different methods.

\begin{conjecture}
We conjecture that one could obtain $\OL_\infty (\A) \le 2$ for every nuclear
$C^*$-algebra  $\A$.
\end{conjecture}

\end{document}